\documentclass[12pt,authoryear]{article}

\usepackage{amsmath,amssymb,amsthm,latexsym,geometry,color,comment,ulem}

\usepackage{mathtools}

\newcommand{\mathsout}[1]
{\bgroup\mathchoice
  {\sbox0{$\displaystyle{#1}$}%
    \usebox0\hspace{-\wd0}%
    \rule[0.5\ht0-0.5\dp0-.5pt]{\wd0}{1pt}}%
  {\sbox0{$\textstyle{#1}$}%
    \usebox0\hspace{-\wd0}%
    \rule[0.5\ht0-0.5\dp0-.5pt]{\wd0}{1pt}}%
  {\sbox0{$\scriptstyle{#1}$}%
    \usebox0\hspace{-\wd0}%
    \rule[0.5\ht0-0.5\dp0-.5pt]{\wd0}{1pt}}%
  {\sbox0{$\scriptscriptstyle{#1}$}%
    \usebox0\hspace{-\wd0}%
    \rule[0.5\ht0-0.5\dp0-.5pt]{\wd0}{1pt}}%
\egroup}
\usepackage[english]{babel}

\newtheorem{Def}{[Def]}[section]
\newtheorem{Prop}[Def]{Proposition}
\newtheorem{Lem}[Def]{Lemma}
\newtheorem{Thm}[Def]{Theorem}
\newtheorem{Cor}[Def]{Cororally}
\newtheorem{Remark}[Def]{Remark}

\newcommand{\Erase}[1]{\if0{#1}\fi}

\begin{document}




\title{A Riemann-Roch theorem on a weighted infinite graph}

\author{Atsushi Atsuji$^\dagger$  \and
Hiroshi Kaneko$^\ddagger$ \and \\
$^\dagger$Department of Mathematics, Keio University, \\
3-14-1 Hiyoshi, Kohoku-ku, Yokohama, Kanagawa 223-8522, Japan \and 
$^\ddagger$Department of Mathematics, Tokyo University of Science, \\ 1-3 Kagurazaka, Shinjuku, Tokyo, 
162-8601, Japan}

\date{}

 


\maketitle

\begin{abstract}
A Riemann-Roch theorem on graph was initiated by M. Baker and S. Norine.
 In their article \cite{BN}, a Riemann-Roch theorem on a finite graph with uniform vertex-weight and 
uniform edge-weight was established and it was suggested a 
Riemann-Roch theorem on an infinite graph was feasible.
In this article, we take an edge-weighted infinite graph and focus on the importance of the spectral gaps of the Laplace operators
defined on its finite subgraphs naturally  given by  $\mathbb Q$-valued positive weights on the edges. 
We build a potential theoretic scheme for
 proof of a Riemann-Roch theorem on the edge-weighted infinite graph.

%
%

\end{abstract}


%


 \section{Introduction}
  A  Riemann-Roch theorem on a connected finite graph was initiated 
 by  M. Baker and S. Norine in \cite{BN}, 
 where connected graph with finite vertices was investigated and 
unit weight was given on each edge of the graph. Since a counterpart of the lowest
exponents of the complex variable in the Laurent series was directly highlighted 
in the Riemann-Roch theorem on graph, research on its relationships with tropical geometry 
was undertaken before other non-tropical 
complex analytical developments were implemented on the graphs.
In fact, A. Gathmann and M. Kerber showed Riemann-Roch theorem on
metric graphs (i.e. graphs where edge lengths are not required to be rational numbers) and
then proceeded 
to
the graphs on tropical curves (i.e. graphs with possibly unbounded edges) in \cite{GK}.
 \medskip
 
In \cite{BS}, M. Baker and F. Shokrieh revealed a close relationships between chip-firing games and potential theory on graphs, 
by characterizing reduced divisors on graphs as the solution to an energy minimization problem.
In their article, an algorithm to find random spanning trees was proposed
where its running times, as well as other algorithms, were analyzed by using potential theoretic ideas.
A new proof of Kirchhoff's matrix-tree theorem was also covered in their article.
R. James and R.  Miranda showed an approach to Riemann-Roch theorem
on an edge-weighted graph by proposing an alternative notion of dimension for divisors
where unit weight is assigned for each vertex.
Recently, S. Backman introduced  in \cite{B}, the notions of edge pivot, cut reversal 
and cycle reversal 
and found a procedure to
produce an acyclic partial orientation as an alternative proof of the Riemann-Roch theorem for a finite graph,
where validity of two stronger hypotheses than (RR1) and (RR2)  in \cite{BN} was given and their importance
was highlighted. On the other hand, Riemann-Roch theorem on an infinite graph has been suggested by M. Baker and S. Norine in \cite{BN}.

  \medskip
 
The objective of this article is to demonstrate a Riemann-Roch theorem on an edge-weighted infinite graph 
where weights on vertices are determined by the sum of the $\mathbb Q$-valued positive weights given on  its adjacent edges. Since the edge-weights give the Laplace operator which is also regarded  as a generator of a Markov process,
the proof is performed by a potential theoretic method including the Dirichlet space theory.
We show a sufficient condition for the spectral gap in terms of the transition probability
of the  Markov process and  prove a Riemann-Roch theorem under certain probabilisitic 
assumptions.

   \medskip
   
In the second section, we reinvestigate the method in \cite{BN} for Riemann-Roch theorem
to accommodate procedures in  \cite{BN}  with $\mathbb Q$-valued edge-weighted finite graph by rethinking 
\textcolor{black}{the canonical divisor in }
the existing Riemann-Roch theorem in \cite{BN}. 
\textcolor{black}{For the sake of self-contained article presentation, 
we attentively look at several assertions updated for a weighted finite graph and pay attention to geometric aspects and characteristics of the weighted graph. As for 
various important relationships between networks with weighted edges and observations in other research subjects, see \cite{K} and \cite{CP}.}
In the third section, \textcolor{black}{since we can find a potential theoretic aspect in the proof of the Riemann-Roch theorem on weighted finite graph, we take validity of Poincar\'e inequality and Riemann-Roch theorem on any edge-weighted finite graph into account and use potential theoretical techniques to approach Poincar\'e inequality on edge-weighed infinite graph.} \textcolor{black}{We focus also on probabilistic significance in 
a decay of edge-weights responding to increase in its distance from a fixed reference point in graph metric basis and
establish a Poincar\'e inequality on the edge-weighted infinite graph by taking an exhaustive sequence of finite subgraphs.} The aim of the section is to validate a spectral gap on the basis of the method as in \cite{Wa} for widely accepted discussions
on  $L^2$-boundedness of $0$-order resolvent
of a probabilistically natural Laplacian associated with the edge-weights.
In the final section, 
we finalize the proof of our Riemann-Roch theorem on the infinite graph
with $\mathbb Q$-valued positive weights on edges.
\textcolor{black}{In the proof, we note that the spectral gap yields
the uniform norm boundedness of a family of functions
on the infinite graph (each function in the family 
represents an accumulation of chip-firings on a subgraph) to obtain a characteristic value of a given
divisor.
Along with a counterpart of the Euler characteristic for the weighted infinite graph,
the characteristic
values of divisors are involved directly in our Riemann-Roch theorem.}

\vspace{1em}
\section{Riemann-Roch theorem on a weighted finite graph}

Let $G=(V_G,E_G)$ be a connected graph consisting of a finite set $V_G$ of vertices and a finite set $E_G$ of edges
 without loops. To be more precise, $E_G$ is given as a subset of $(V_G\times V_G\setminus \{\{x,x\} \mid x \in V_G\})/\sim$ \textcolor{black}{, where
 $\{x,y\} \sim \{y,x\}\in E_G$ for any pair $x,y\in V_G$. Accordingly, $\{x,y\}$ and $\{y,x\}$ are regarded as an identical element in $E_G$.} We assume that a $ \mathbb Q$-valued positive
weight $C_{x,y}$ is given at every edge $\{x,y\}\in E_G$ \textcolor{black}{with the symmetricity $C_{x,y}=C_{y,x}$ for any $\{x,y\} \in E_G$, according to the identification $\{x,y\}=\{y,x\}$.}
At each point $x$ in $V_G$, the set $N(x)$ consisting of its neighbor is defined by $N(x)=\{y\in V_G\mid \{x,y\} \in E_G\}$ of $x$.
For any $\mathbb Z$-valued function $f$ on $ V_G$,  \textcolor{black}{the Laplacian} $\Delta f$ \textcolor{black}{ of the function} is defined by $\Delta f(x)=\sum_{y \in N(x)} C_{x,y}(f(x)-f(y))$ as \textcolor{black}{ a function} on $V_G$.
\textcolor{black}{When $G$ is regarded as the electrical circuit given by installing a resistor with the resistance $1/C_{x,y}$ at every edge $\{x,y\} \in E_G$,}
the value $f(x)$ is viewed as the voltage level at $x \in V_G$. \textcolor{black}{We can focus on}
another function $i$ on $V_G$ 
defined by $i(x)=\min \{\vert \Delta f(x)\vert  \mid f: V_G\to \mathbb Z \mbox{ satisfying $f(x)=0$ and }\Delta f(x) \not=0\}$, which is viewed as
the minimum positive absolute value 
of \textcolor{black}{ feasible} current flows at \textcolor{black}{the grounded vertex }$x$ \textcolor{black}{ given} by integer valued voltage $f(y)$ with $y \in N(x)$.
In what follows\textcolor{black}{, similarly to M. Baker and S. Norine's article \cite{BN}, }
 $f$ \textcolor{black}{ persistently }stands for $\mathbb Z$-valued function on $V_G$.
\medskip

We use notation $\ell$ for another $\mathbb Z$-valued function
to assign an integer multiple of $i(x)$ at each $x \in V_G$. 
A divisor on the graph $G$ is given by $D=\sum_{x \in V_G} \ell(x)i(x)1_{\{x\}}$
and its degree ${\rm deg}(D)$ is defined by ${\rm deg}(D)=\sum_{x \in V_G} \ell(x)i(x)$.
\textcolor{black}{ A divisor $D=\sum_{x \in V_G} \ell(x)i(x)1_{\{x\}}$ is said to be effective if $\ell(x)\geq 0$ for all $x \in V_G$.
Since $f$ is $\mathbb Z$-valued, $\Delta f(x)$ is given as an integer multiple of $i(x)$ at each $x \in V_G$ and the Laplacian}
$\Delta f$
will be identified with the divisor $\sum_{x \in V_G} \Delta f(x)1_{\{x\}}$. \textcolor{black}{This identification makes it possible to
add the Laplacian $\Delta f$ to any divisor.}
\medskip

\textcolor{black}{Along with this addition, we need two materials characterizing the graph $G$
for a Riemann-Roch theorem on the weighted graph: 
one is the positive real value 
$\min \{\vert \sum_{x\in V_G}\ell(x)i(x)\vert\in (0,\infty)\mid \ell: V_G\to \mathbb Z\}$, which is denoted 
by $i_{(G,C)}$; and the other is}
the canonical  divisor $K_G$ on the weighted graph $G$, \textcolor{black}{which} is given by $K_G=\sum_{x \in V_G} \{\sum_{y \in N(x)}C_{x,y}-2i(x)\}1_{\{x\}}$.
\textcolor{black}{ Similar to \cite{BN}, for our goal, we take advantage of 
total orders on $V_G$ and denote the family of total orders on $V_G$ by $\mathcal O$.  For each $O \in \mathcal O$, }
we introduce the divisor $\nu_O$ given by

$$\sum_{x \in V_G} \nu_O(x)1_{\{x\}},$$

\noindent
with $\nu_O(x)= \sum_{y\in N(x), y<_O x}C_{x,y}-i(x) \mbox{ for }x \in V_G$,
and its reversed total order $\overline O$ 
 is
 defined by $x<_{\overline O}y$ for any $x,y \in V_G$ satisfying $y<_{O}x$.
 
 \medskip
 \textcolor{black}{The importance of such divisors is found in its close relationship
 with the canonical  divisor $K_G$ as shown in the first assertion of the following lemma, similar to  \cite{BN}.}

\begin{Lem}\label{deg} For any divisor $\nu_{\overline O}$ determined by $\overline O \in \mathcal O$ and  any divisor $\Delta f$ 
determined by $\mathbb Z$-valued function $f$, one sees the followings:
\begin{itemize}
\item[\rm{(i)}] $\nu_{\overline O}=K_G-\nu_{O}$,
 \item[\rm{(ii)}] ${\rm deg}(\Delta f)=0$.
 \end{itemize}
\end{Lem}

{\it Proof.}\quad  (i) The assertion is clear from the identity $\nu_{\overline O} + \nu_{O}=K_G$. (ii)
The identity follows from the equivalence between $y\in N(x)$ and $x\in N(y)$
\textcolor{black}{ and the trivial identity $C_{x,y}(f(x)-f(y))=-C_{y,x}(f(y)-f(x))$.}\qed
\medskip

If a pair of  divisors $D$ and $D'$ 
satisfies
$D'=D+\Delta f $ for some $\mathbb Z$-valued function $f$, $D$ and $D'$ are called equivalent
and the relationship is denoted by $D \sim D'$. As in \cite{BN}, we introduce \textcolor{black}{ the linear system }
$\vert D\vert =\{D' \mid D'\mbox{ is effective and equivalent with }$ $ D\}$ for any divisor $D$ and the following condition on the graph $G=(V_G,E_G)$:
\medskip

\noindent
{\bf (RR)} For each divisor $D$, there exists an $O \in \mathcal O$ such that \textcolor{black}{exactly one of two linear systems $\vert D \vert$ 
and 
$\vert \nu_O-D \vert$} is empty.

\medskip
\noindent
The family of total orders satisfying this condition \textcolor{black}{ with respect to } $D$ will be denoted by $\mathcal O_D$. 
\medskip

We introduce 
$\mathfrak e_{(G,C)}=\sum_{x\in V_G} i(x)-\sum_{\{x,y\} \in E_G}C_{x,y}$ as a counterpart of $1-g$ with the genus $g$ of $G$  in M. Baker and S. Norine's article  \cite{BN}. \textcolor{black}{ In fact, we see that}

$$\textcolor{black}{{\rm deg}(K_G)=\sum_{x \in V_G} \{\sum_{y \in N(x)}C_{x,y}-2i(x)\}=-2 \mathfrak e_{(G,C)}}.$$
\medskip

\begin{Lem}\label{notequivalent} For each  $O\in \mathcal O$, $\nu_O \in \mathcal N$, where $\mathcal N$ stands for the family of divisors of degree $-\mathfrak e_{(G,C)}$ admitting only
non-effective equivalent divisors.
\end{Lem}
\medskip

{\it Proof.}\quad It is easy to see that
 $\mbox{deg}(\nu_O)=-\mathfrak e_{(G,C)}$ and $\nu_O$ is not effective as in \cite{BN}. For any divisor $D$ \textcolor{black}{given by}
$D = \nu_O -\Delta f$ \textcolor{black}{ with } some non-constant $\mathbb Z$-valued function $f$, we take $A_f=\{x \in V_G \mid f(x) = \max_{y \in V_G}f(y)\}$
and the minimal element $z$ in $A_f$ 
with respect to the 
total
order $O$. Since $y <_O z$ implies $f(y) \leq  f(z)-1$ and $f(y) \leq  f(z)$ for all $y \in V_G$,
it turns out that

\begin{align*}
D(z) &= \Big(\sum_{y \in N(z),y <_O z}C_{z,y}-i(z)\Big) - \Big(\sum_{y\in N(z)}C_{z,y}(f(z)-f(y))\Big)\\
&=-i(z)+\sum_{y\in N(z),y >_O z}C_{z,y}(f(y)-f(z))+\sum_{y\in N(z),y <_O z}C_{z,y}(f(y)-f(z)+1)
\\
&\leq-i(z).
\end{align*}

\noindent
This implies that $\nu_O$  can not be equivalent to any effective divisor. 
 \qed
 
 \medskip
 
 For any divisor $D$ and non-negative integer $k$, we introduce  $E_k(D) =\{\mbox{ effective }\linebreak\mbox{divisors }
E\mid \mbox{deg }(E)=k\,i_{(G ,C)} \mbox{ satisfying }\vert D-E\vert\not=\emptyset\}$. We note that
$E_0(D)$ is either the empty set or the set consisting only of the \textcolor{black}{ zero} divisor \textcolor{black}{ vanishing on $V_G$}.
We also define a $\{-
i_{(G ,C)},0, i_{(G,C)},2 i_{(G,C)},\dots \}$-valued function $r$  on the set consisting of 
all divisors by
$$
r(D)=\begin{cases} 
-i_{(G,C)},  \qquad& \mbox{ if $E_0(D)=\emptyset$},  \\
\max\{k\,i_{(G,C)} &  \mid  E_k(D)  \mbox{ consists of all  effective divisors of degree }k\,i_{(G,C)}\}, \\
& \mbox{ otherwise},
\end{cases}
$$

\noindent
similar to \cite{BN}. 
 \medskip

\textcolor{black}{ Our objective is to establish the following identity as a natural counterpart of the Riemann-Roch theorem on a weighted 
graph:}

$$\textcolor{black}{r(D)-r(K_G-D) = {\rm deg}(D) + \mathfrak e_{(G,C)}.}$$

 \medskip
 
Let us briefly look back at the procedures for the proof of Riemann-Roch theorem presented in the
article based on the notions in our settings \textcolor{black}{involving edge-weight}.

 \medskip

 For a divisor $D=\sum_{x \in V_G} \ell(x)i(x)1_{\{x\}}$, the subset $\{x \in V_G\mid \ell(x)\not=0\}$ is called a support of $D$ and denoted by {\rm supp}$[D]$. It is straightforward that any  divisor $D$ is uniquely decomposed as $D=E-E'$ with effective divisors $E,E'$
 without intersection of their supports. In what follows, the degree 
 of $E$ in the decomposition will be denoted by ${\rm deg}^+(D)$
 and the degree of $E'$ in the decomposition will be denoted by ${\rm deg}^-(D)$.
  \medskip
 
 \begin{Prop}\label{representationforr} Condition {\bf (RR)} implies that 

$$r(D)=\Big( \min_{D'\sim D,  O \in \mathcal O}\mbox{\rm deg}^{+}(D' - \nu_O)\Big)-i_{(G,C)}$$

\noindent
for any divisor $D$.

\end{Prop}

{\it Proof.}\quad  If $\vert D - E\vert=\emptyset$ for some effective divisor $E$, {\bf (RR)} implies that there exists a total order 
$O \in \mathcal O_{D-E}$ 
such that $\nu_O -(D-E) \sim E'$ equivalent to $\nu_O -D \sim E'-E$ with some effective divisor $E'$. This shows that $D' - \nu_O \sim E-E'$ 
with the effective divisors $E,E'$ independent of the choice of $D'\sim D$.
 \medskip
 
Conversely, if $D' - \nu_O \sim E-E'$ for some effective divisors $E,E'$,
$O \in \mathcal O$ and \textcolor{black}{$D'$ satisfying} $D'\sim D$, then $D-E\sim \nu_O- E'$. 
However, Lemma \ref{notequivalent} implies that $\nu_O - E'$ \textcolor{black}{ is }not 
equivalent to any effective divisor. This shows that $\vert D -E\vert=\emptyset$.
 \medskip

Originally, 
\textcolor{black}{ by taking $I_{(G,C)}=\{k\, i_{(G,C)} \mid k\in \{0,1,2,\dots\}\}$,}
the value $r(D)$ is characterized by
$r(D)<s$ if and only if there exists some effective divisor $E$ with ${\rm deg}(E)=s$ such that $\vert D - E\vert=\emptyset$,
\textcolor{black}{ where $s \in I_{(G,C)}$.}
Another charactarization can now be admitted
$r(D)<s$ if and only if $s\geq \mbox{\rm deg}^{+}(D' - \nu_O)$ for some $D'\sim D$ and $O \in \mathcal O$,
\textcolor{black}{ by using non-negative integer multiple $s$ of $i_{(G,C)}$.}
The assertion has been proved. \qed
\medskip

This proposition shows that, for a divisor $D$ on $V_G$, its equivalent  divisor $D'$ and a total order $O$ on $V_G$ are taken so that 
the minimum in the proposition is attained.  We introduce the involution $\overline D'$ of $D'$ by $\overline D'= K_G - D'$
and \textcolor{black}{ then see }that $\overline D' - \nu_{\overline O}= K_G - D' - (K_G - \nu_O)=\nu_O -D' $, which gives an alternative
expression $\overline D' - \nu_{\overline O}\sim E'-E$, the decomposition into the sum of effective divisors \textcolor{black}{ as }in the 
proof of the last proposition.
\medskip

\begin{Cor}\label{minimaldegree} Condition {\bf (RR)} implies

\begin{itemize}

\item[\rm{(i)}] for any divisor $D$, 
$$r(\overline D)=\Big( \min_{D'\sim D,  O \in \mathcal O}\mbox{\rm deg}^{+}(\overline D'-\nu_{\overline O})\Big)-i_{G 
,C)},$$
where $\overline D = K_G - D$, the involution of $D$,
\item[\rm{(ii)}] $\min\{ {\rm deg}(E) \mid D' - \nu_O \sim E-E' \mbox{ with some effective divisors }E,E'\}$ and $\min\{ \linebreak {\rm deg}(E') \mid \overline D' - \nu_{\overline O} \sim E'-E\mbox{ with some effective divisors }E,E'\}$ are \textcolor{black}{ both }attained \textcolor{black}{by the effective divisors $E=E_0,E'=E'_0$}, which are characterized by 
a unique decomposition 
$D' - \nu_O \sim E_0-E'_0$ into the sum of effective divisors of minimal degrees, equivalently given as $\overline D' - \nu_{\overline O} \sim E'_0-E_0$,
\item[\rm{(iii)}] for any divisor $D$ with  ${\rm deg}(D)=-\mathfrak e_{(G,C)}$, 
$D \in \mathcal N$ if and only if $K_G-D \in \mathcal N$.
\end{itemize}
\end{Cor}

\medskip

{\it Proof.}\quad (i) The identity follows from the equivalence between $D'\sim D$ and $\overline D'\sim \overline D$
and from the equivalence between { the }choice of $O$ and \textcolor{black}{the choice }of $\overline O$.
(ii) Since ${\rm deg}(E)=\mathfrak e_{(G,C)} +{\rm deg}(D) +{\rm deg}(E')$, the assertion is 
straightforward and justified. (iii) The assertion follows from ${\rm deg}(K_G-D)=-\mathfrak e_{(G,C)}$ and the equivalence 
between $D\notin \mathcal N$ and $K_G-D =\overline D\notin \mathcal N$, i.e., 
the equivalence between $ r(D')=r(D) \geq 0$ and $r(\overline D')=r(\overline D) \geq 0$.
\qed

\medskip

\begin{Remark} Condition {\rm (RR2)} in  {\rm \cite{BN}} is affirmatively asserted in {\rm (iii)} in the corollary.
\end{Remark}
\medskip

\textcolor{black}{We introduce the graph metric $d$ which is defined by}
$$\textcolor{black}{d(v,w) = \begin{cases}\min\{j \mid \{w_0,w_1\}, \dots, \{w_{j-1},w_j\} \in E_G, v=w_0 \mbox{ and }w=w_j\} & \mbox{ if }v\not=w,\\
0 & \mbox{ if }v =w.
\end{cases}}$$

\medskip
\textcolor{black}{We take a fixed base (reference) vertex $v_0 \in V_G$, then the distance between $v$ and $v_0$ is given by
$d(v_0,v)$, which will be denoted by $d_G(v)$.} 
\medskip

\textcolor{black}{We introduce $d_G= \max\{d(v) \mid v \in V_G\}$, and} for any given divisor $D=\sum_{x \in V_G} \ell(x)i(x)1_{\{x\}}$ on $V_G$, we introduce a $d_G-1$-dimensional vector

$$\mathbf V_1(D)= \Big(\sum_{z \in S_{d_G}, \ell(z) <0}\ell(z)i(z),\sum_{z \in S_{d_G-1}, \ell(z) <0}\ell(z)i(z),\dots, \sum_{z \in S_{1},\ell(z) <0}\ell(z)i(z)\Big)$$

\noindent
and a $d_G$-dimensional vector

$$\mathbf V_2(D)= \Big(\sum_{z \in S_0}\ell(z)i(z),\sum_{z \in S_{1}}\ell(z)i(z),\dots, \sum_{z \in S_{d_G}}\ell(z)i(z)\Big),$$

\medskip

\noindent
\textcolor{black}{where $S_j=\{v \in V_G \mid d_G(v)=j\}$ for any non-negative integer $j$.}

\medskip

If a divisor $D=\sum_{x \in V_G} \ell(x)i(x)1_{\{x\}}$ satisfies the following two conditions 

\begin{itemize}
\item[(P1)] $\ell(z)i(z) \geq 0 \mbox{ for all }z\in V_G \setminus\{v_0\}$,
\item[(P2)] for every non-empty set $A\subset V_G \setminus\{v_0\}$, there exists a vertex
$x\in A$ such that $\ell(x)i(x) < \mbox{{\rm outdeg}}_A(x)=\sum_{y \in N(x)\cap A^c}C_{x,y}$,
\end{itemize}

\noindent
with respect to the base vertex 
$v_0$,  then $D$ is said to be  $v_0$-reduced.
\medskip

\begin{Prop} For any divisor $D$, there exists
a unique $v_0$-reduced divisor $D_0$ such that $D_0 \sim D.$
\end{Prop}

{\it Proof.}\quad \textcolor{black}{By starting with the originally given divisor $D$, we first take the subfamily $\mathcal D'$ of its equivalent divisors given
by}
\begin{align*}
\mathcal D'=\{D' \mid & D' \mbox{ attains the maximum }\max_{D'\sim D}\mathbf V_1(D')\mbox{ in the sense of the lexicographical} \\
& \mbox{ order}\},
\end{align*}

\noindent
and in the next, we take the subfamily $\mathcal D''$ of $\mathcal D'$ given by
\begin{align*}
\mathcal D''=\{D'' \in \mathcal D' \mid & D'' \mbox{ attains the maximum }\max_{D''\sim D}\mathbf V_2(D'')\mbox{ in the sense of the } \\
& \mbox{  lexicographical order}\}.
\end{align*}

\noindent
We see that any divisor in $\mathcal D''$ attains both maxima and satisfies the conditions (P1) and (P2) with respect to 
the base (reference) vertex 
$v_0$ as in the proof of Proposition 3.1 in \cite{BN}.
\medskip

For the uniqueness of the $v_0$-reduced divisor, we assume that there exist two distinct $v_0$-reduced divisors $D=\sum_{x \in V_G} \ell(x)i(x)
1_{\{x\}}$, $D'=\sum_{x \in V_G} \ell'(x)i(x)1_{\{x\}}$ such that $D' = D - \Delta f$ for some non-constant $\mathbb Z$-valued function $f$. Without losing these general settings, we may assume that 
$\max_{y \in V_G} f(y)$ $> f(v_0)$
as in the proof of uniqueness in the proposition in \cite{BN}, from which we can derive that
$v_0 \notin A_f=\{x\in V_G \mid f(x)=\max_{y \in V_G} f(y)\}$ and 

\textcolor{black}{
$$
0 \leq \ell'(x)i(x) = \ell(x)i(x) -\sum_{\{x,y\} \in E_G}C_{x,y}(f(x) -f(y))\leq  \ell(x)i(x)-{\rm outdeg}_{A_f}(x),
$$
}
\noindent
for any $x \in A_f$.  This contradicts the condition (P2) originally imposed on $D$. 
\qed 
 
 \medskip
 The following theorem shows that (RR1) in \cite{BN} still holds valid even for our weighted graph.
  \medskip

\begin{Thm}\label{nonemptyOD} For any divisor $D$, $\mathcal O_D\not=\emptyset$, namely Condition {\bf (RR)} is satisfied.
\end{Thm}
\medskip

{\it Proof.}\quad A similar procedure to the construction of a total order in  the proof of Theorem 3.3 in \cite{BN} 
works to deduce $\mathcal O_D\not=\emptyset$. \qed
\medskip

\begin{Remark} For any divisor $D$ with ${\rm deg}(D)=-\mathfrak e_{(G,C)}$, $D \in \mathcal N$ if and only if $D \sim
\nu_O$ for some $O \in \mathcal O$.
In fact, thanks to Lemma \ref{notequivalent},  it suffices to show that $D \in \mathcal N$ implies the existence of $O \in \mathcal O$ satisfying $D \sim \nu_O$. This implication is verified by combining $\mathcal O_D\not=\emptyset$
with the fact that any $D \in \mathcal N$ does not admit any equivalent effective divisors. \end{Remark}

\begin{Remark}  Our weight function $i(x)$ defined on the set $V_G$ of vertices 
may differ from the uniform vertex-weight assigned 
in  {\rm\cite{BN}}. On the other hand, we see that
$C_{x,y}$ admits representation as $C_{x,y}=\dfrac{n_{x,y}}{m_G}$ 
with positive integers $n_{x,y}$ and $m_G$, where $m$ can be taken independently of 
edges $\{x,y\}\in E_G$.
Accordingly, we can
install $n_{x,y}$-ple edges with a weight between $x$ and $y$ in  the network discussed in 
{\rm\cite{BN}} so that essentially the same network is obtained up to the difference in minimum edge weights,
namely, such a difference between $i_{(G,C)}$ and $1/m_G$.
This shows us that the procedure in {\rm\cite{BN}} for finding the $v_0$-reduced divisor is still valid. 
This is because the image $\{\Delta f \mid f: V_G \to \mathbb Z\}$ of our Laplacian is linearly isomorphic to the one in  {\rm\cite{BN}}.
 \end{Remark}
\medskip

\begin{Thm}{\bf (Riemann-Roch theorem on a weighted finite graph).}  For any divisor $D$, 

$$r(D)-r(K_G-D) = {\rm deg}(D) + \mathfrak e_{(G,C)}.$$
\medskip

\end{Thm}

{\it Proof.}\quad Theorem \ref{nonemptyOD} shows that $\mathcal O_D\not=\emptyset$ for any divisor $D$. Therefore, one
can take a divisor 
\textcolor{black}{$D'$ satisfying }$D'\sim D$ and a total order $O \in \mathcal O_{D-E_0}$ for any effective divisor $E_0$ satisfying
$\vert D-E_0\vert=\emptyset$ \textcolor{black}{and} the minimal degree condition 
${\rm deg}(E_0)=\min\{ {\rm deg}(E)\mid \vert D-E\vert=\emptyset, E \mbox{ is effective}\}$.
By the discussion in the proof of Corollary \ref{minimaldegree}, we have the decomposition
$D' -\nu_O =E_0-E'_0$ with effective divisors $E_0, E'_0$ of the minimal degree.
Accordingly, we see that

\begin{align*}
{\rm deg}^+(D'-\nu_O) -{\rm deg}^+(\overline D'-\nu_{\overline O})&={\rm deg}^+(D'-\nu_O) -{\rm deg}^+((K_G -D')-(K_G -\nu_O)) \\
& = {\rm deg}^+(D'- \nu_O) - {\rm deg}^+(\nu_O -D')\\
&={\rm deg}(E_0)-{\rm deg}(E'_0)\\
&= {\rm deg}(D' -\nu_O)\\
&= {\rm deg}(D) + \mathfrak e_{(G,C)}.
\end{align*}
\medskip

\noindent
Since the left-hand side is equal to $r(D)-r(K_G-D)$ due to Proposition \ref{representationforr},
the identity in the assertion is derived.\qed
\medskip
 
 We close this section with the following fundamental properties of $r(D)$ which are utilized 
 later.
  \medskip

\begin{Lem}\label{fundamental}
{\rm (i)}\quad If $D'$ is effective, then $ r(D)+{\rm deg}(D')\geq r(D+D')$,
\smallskip

\noindent
{\rm (ii)}\quad If $-D''$ is effective, then $r(D)+{\rm deg}(D'')\leq r(D+D'')$.
\end{Lem}

{\it Proof.}\quad (i) Since there exists an effective divisor $E$ with ${\rm deg}(E)> r(D)$ satisfying
$\vert D-E\vert =\emptyset$, it turns out that the effective divisor $D'+E$ satisfies ${\rm deg}(E+D')> r(D)+{\rm deg}(D')$ and $\vert (D+D')-(E+D')\vert =\emptyset$.
Accordingly, $r(D+D')$ does not exceed $r(D)+{\rm deg}(D')$. \\
\smallskip

\noindent
(ii) It suffices to show that $r(D)\leq r(D+D'')-{\rm deg}(D'')$. This follows from (i) by regarding $D+D''$
and $-D''$
 as $D$ and $D'$ respectively in the identity of (i).\qed
\medskip


 \section{Poincar\'e inequality for spectral gap on an infinite graph}
  \medskip
   We shift our attention from finite graph{s} to {an} infinite graph.
  Throughout this section, we consider a connected infinite graph $G=(V_G,E_G)$ consisting 
  of countably infinite sets $V_G$ and $E_G$ of vertices and edges, respectively. We assume that 
  the graph $G=(V_G,E_G)$ does not admit any loops and 
\textcolor{black}{\Erase{$G$} is locally finite, namely 
$N(x)=\{y \in V_G \mid \{x,y\} \in E_G\}$ is a finite set
for every $x \in V_G$. 
Let $C_{x,y}$ denote $ \mathbb Q$-valued positive
weight given at every edge $\{x,y\}\in E_G$ satisfying that $C_{x,y}=C_{y,x}$ for every $x,y\in V_G$, according to 
the identification $\{x,y\}=\{y,x\}$. Define a measure $m$ on $V_G$ by
\[
m(A)= \sum_{x\in A} \sum_{y\in N(x)}C_{x,y},
\]
for $A \subset V_G$. 
We also assume finiteness of total volume on the graph $G$, namely, $m(V_G)<\infty$.
}
\medskip

{\color{black}{
Here, we describe the basic idea of obtaining  
our Riemann-Roch theorem on an infinite graph from the facts obtained in the case of finite graphs. 
We take the graph metric $d$ \textcolor{black}{that is }exactly \textcolor{black}{the} same \textcolor{black}{in} the previous section
 and fix a base \textcolor{black}{(reference)} vertex $v_0$, define $V_n=\{x \in V_G\mid d(v_0,x) \leq n\}$ and $E_n=\{\{x,y\} \mid x,y \in V_n\}$ to determine the subgraph $G_n=(V_n,E_n)$ and introduce the Laplacian $\Delta_n f(x)=\sum_{y\in N(x)\cap V_n}C_{x,y}(f(x)-f(y))$ on $V_n$.

\medskip

As will be seen later, a divisor $D$ can be given globally on $G$ and it admits
a sequence $\{(D)_n\}$ of divisors, $n$-th divisor of  which is a restriction of $D$ to $G_n$.  We already have the Riemann-Roch theorem on $G_n$ in Section 2. \textcolor{black}{Thus}, 
\textcolor{black}{ by the equivalence $D' \overset {n}{ \sim} D''$ between divisors $D'$ and $D''$ on $V_n$ defined by $D''=D' + \Delta_n f$ for some $\mathbb Z$-valued function $f$ on $V_n$
and with the same notations in the last section,
we have $r_n((D)_n)$ on $G_n$ satisfying }
 \begin{equation*}
r_{n}((D)_n)
=\Big( \min_{D'\overset {n}{ \sim} (D)_n, 
O_{n} \in \mathcal O_{n}}
{\rm deg}^{+}(D' - \nu_{O_{n}})\Big)-i_{(G_{n},C_{n})},
\end{equation*}
\medskip
 \noindent
\textcolor{black}{where $\mathcal O_{n}$ stands for the set consisting of all total orders on $V_n$}. 
To be more precise,  %
%
bearing the difference between $\Delta_n f(x)$ and $\Delta f(x)=\sum_{y\in
N(x)}C_{x,y}(f(x)-f(y))$
for $x$ with $d(v_0,x)=n$ in our mind, we notice that
$\sum_{x \in V_{n}} \ell(x)i(x)1_{\{x\}}$ can not always be regarded
as any divisor on $G_n$, and in contrast
$\sum_{x \in V_{n-1}} \ell(x)i(x) 1_{\{x\}}$ is regarded as a divisor
vanishing outside $V_{n-1}$
as detailed in 
Remark \ref{restrictionofdivisor} below.
We take the latter  vanishing outside $V_{n-1}$ as the restriction $(D)_n$
to $G_n$, when a divisor $D=\sum_{x \in V_{G}} \ell(x)i(x) 1_{\{x\}}$ is
globally given.

 \begin{Remark}\label{restrictionofdivisor}\quad For any positive integer $n$ and $x \in V_n$, not only $i(x)=\min \{\vert \Delta f(x) \vert \mid f: V_G\to \mathbb Z$ with $\Delta f(x)\not=0\}$
 but $i_{n}(x)=\min \{\vert \Delta_n f(x) \vert \mid f: V_n\to \mathbb Z $ with $\Delta_n f(x)\not=0\}$ is defined. 
We note that 
\textcolor{black}{$i(x)$ does not always coincide with $i_n(x)$ for $x$ with $d(v_0,x)=n$\textcolor{black}{,} and  that
$D'=\sum_{x \in V_{n}} \ell(x)i(x)1_{\{x\}}$ does not always admit \textcolor{black}{the} expression $\sum_{x \in V_{n}} \ell(x)i_n(x)1_{\{x\}}$. However,
 for any pair of positive integer\textcolor{black}{s} $j,n$ with $j\geq n$, 
 any divisor $D=\sum_{x\in V_{n-1}}\ell(x)i(x)1_{\{x\}}$ vanishing outside $V_{n-1}$ is regarded as a divisor on $G_j$.
}
%
 \end{Remark}
 \medskip
 
\textcolor{black}{An} integer-valued function $f_n$ on $V_{n}$\textcolor{black}{,} such that $D'= (D)_n + \Delta_{n}f_n$ attains the minimum
 $\min_{D'\overset {n}{ \sim} (D)_n, 
O_{n} \in \mathcal O_{n}}
{\rm deg}^{+}(D' - \nu_{O_{n}})$\textcolor{black}{,} is called \textcolor{black}{a} minimizer for $r_{n}((D)_n)$. 
\Erase{Here\textcolor{black}{,} $\Delta_n$ is the weighted Laplacian on $G_n$\textcolor{black}{,} whose basic properties we will see in Section 3.1 below. } 
Then\textcolor{black}{,} one of \textcolor{black}{the problems to be considered} here is the \textcolor{black}{global} existence of such a function\Erase{ globally}, namely, the solution $f_n$ to the equation $\Delta_{n}f=D'-(D)_n$ and extension of the function to larger subgraph.
%
In the case of finite graph, we always have the solution 
\textcolor{black}{to $\Delta_n f=g$ on $G_n$, since the solution can be given by 
the $0$-order resolvent of $\Delta_n$. 
If  we can expect a stability of the correspondence from $D'-(D)_n$ to $f_n$ as $n\to\infty$, 
using $0$-order resolvent as a main tool still in infinite graph, 
} 
%
we\Erase{ will} carry out the following procedure: 
First, we take a minimizer $f_n$ by using resolvent of $\Delta_n$ on $G_n$. Second, by deriving 
the convergence of $\{f_n\}$
as $n\to\infty$ from the stability, we discuss the existence of $\lim_{n\to\infty}{r_n((D)_n)}$, and finally we define $r(D)$ as the limit and complete a proof of
 \textcolor{black}{a} 
Riemann-Roch theorem on infinite graphs. 
\medskip

However, \Erase{actually} these procedures are not \textcolor{black}{actually}  straightforward and we need \textcolor{black}{to carefully discuss} 
 the convergence 
by selecting subsequences $\{(D)_n\}$ and $\{f_n\}$, linked with a choice of subsequence of
$\{O_{n}\}$ with $O_{n} \in \mathcal O_n$. Precise discussion on this problem will be given in Section 4. 
\medskip

The more basic problem here is that the $0$-order resolvent does not always exist on infinite graphs and 
operator theoretical treatment of the resolvent $\Delta_n^{-1}$ for our purpose is unclear within general operator theory. 
The existence of \textcolor{black}{a} {\it spectral gap }
of the Laplacian are known as a powerful condition for the existence of the resolvent (see \cite{Wa} for the general theory of spectral gaps and related functional inequalities) and we can also clear the convergence problem. 
On the other hand, the resolvent can be defined from a natural Markov process on the graphs associated with the weight $C_{x,y}$, which we will give in Section 3.1. To ensure the existence of the spectral gap, we propose a condition of \Erase{our }infinite graphs using a quantity $\rho_n$ on $G_n$ in terms of the Markov process as defined in (\ref{rho_0}) below. In fact, we will show in Section 3.2 that the spectral gap exists if $\limsup_{n\to\infty}\rho_n$ is sufficiently small. Roughly speaking, this condition implies a strong recurrence property of the process and then a graph satisfying the condition looks 
 \textcolor{black}{closely like} 
a finite graph from the view of the Markov process. To show the existence of a spectral gap of our Laplacian under the condition, we introduce a Dirichlet form $\mathcal E$ associated with $\Delta$ in Section 3.1 and show a Poincar\'e type inequality of $\mathcal E$ over some $L^2$ space in Section 3.2. We remark that there are one to one correspondence \Erase{relationships }among the Markov process, the Laplacian and the Dirichlet form. 
}}

%
  
 \subsection{Basic properties of weighted Laplacians} 
In this section, we  use the following notions and notations. 
Let $U=(V_U,E_U)$ be a connected subgraph of $G$, i.e. $U$ is a graph consisting of
set $V_U$ of vertices with $V_U \subset V_G$ and set $E_U$ of edges  with $E_U=\{\{x,y\}\in E_G \mid x,y \in V_U\}$, 
and  $m_U$ be the measure on $V_U$ 
determined by $m_U(A)= \sum_{x\in V_U\cap A} m_U(\{x\})$
\Erase{for any $A\subset U$} with $ m_U(\{x\}) =\sum_{y \in N(x)\cap V_U}C_{x,y}$.
Here and throughout this subsection, we denote $V_G$ by $V$.
\medskip

Recall the definition of the Laplacians on graphs $G$ and $U$;
\begin{align*}
\Delta\phi(x)&=\sum_{y\in N(x)}C_{x,y}(\phi(x)-\phi(y)), \\
\Delta_U u(x)&=\sum_{y\in N(x)\cap V_U}C_{x,y}(u(x)-u(y)),
\end{align*}
\Erase{respectively}\textcolor{black}{for real-valued functions $\phi$ and $u$ on $V$ and $V_U$ respectively.} 

\medskip

Define the operators by
\[
L=\frac1{m}\Delta \mbox{ and } L_U=\frac1{m_U}\Delta_U.
\]
Here and in the sequel, $m$ and $m_U$ stand for the function taking $m(\{x\})$ at every $x$ and  the one taking
$m_U(\{x\})$
at every $x\in U$, respectively.
In what follows, $m(\{x\})$ and $m_U(\{x\})$ will be denoted by $m(x)$ and $m_U(x)$, respectively.

\medskip

Note that $L$ and $L_U$ are symmetric, respectively in the sense that
\[
(L\phi,\psi)_{L^2(m)}=  (\phi, L\psi)_{L^2(m)} \mbox{ and }(L_U u,v)_{L^2(m_U)} =(u, L_U v)_{L^2(m_U)}
\]
for $\phi,\psi \in L^2(m)$ with $L\phi,L\psi \in L^2(m)$ and for $u,v \in L^2(m_U)$ with $Lu,Lv \in L^2(m_U)$.
\medskip

We see that those operators are associated with Dirichlet forms:

\[
{\mathcal E}(f,g)=\sum_{\{x,y\}\in E_G} C_{x,y}(f(x)-f(y))(g(x)-g(y))  
\]
for functions $f,g$ on $V$ with $f,g \in \mathcal F=\{h \in L^2(m) \mid 
{\mathcal E}(h,h)=
\sum_{\{x,y\}\in E_G}C_{x,y}(h(x)-h(y))^2<\infty\}${,} 
and

\[
{\mathcal E}_U(u,v)= \sum_{\{x,y\}\in E_U,\; x,y\in \textcolor{black}{V_U}} C_{x,y}(u(x)-u(y))(v(x)-v(y))
\]
for functions $u,v$ on $U$ with $u,v \in \mathcal F_U=\{w \in L^2(m_U) \mid 
{\mathcal E_U}(w,w)=
\sum_{\{x,y\}\in E_U}C_{x,y}\linebreak (w(x)-w(y))^2<\infty\}$,
in the sense that 
\medskip

\[
(Lf,g)_{L^2(m)}={\mathcal E}(f,g) \mbox{ and }(L_U u,v)_{L^2(m_U)}={\mathcal E}_U(u,v).
\]
\medskip

\noindent
\Erase{as long as both sides make sense.}
\medskip

\Erase{By identifying function $\phi$ on $V$ with its restriction to $U$, we also have}
\textcolor{black}{For a function $\phi$ on $V$} 
\begin{align*}
L\phi(x)&=\frac1{m(x)}\sum_{y\in N(x)}C_{x,y}(\phi(x)-\phi(y)) \\
&=\frac1{m(x)}\sum_{y\in N(x)\cap V_U}C_{x,y}(\phi(x)-\phi(y))+ \frac1{m(x)}\sum_{y\in N(x)\cap V_{U^c}}C_{x,y}(\phi(x)-\phi(y)) \\
&=\frac{m_U(x)}{m(x)}L_U \phi(x) + \frac1{m(x)}\sum_{y\in N(x)\cap V_{U^c}}C_{x,y}(\phi(x)-\phi(y)).
\end{align*}

\noindent
\Erase{From this,}  
\textcolor{black}{By identifying function $\phi$ on $V$ with its restriction to \textcolor{black}{ $V_U$},}  we have  the following relationship between $L$ and $L_U$ {(}which is used frequently {below)}: 
\begin{equation}\label{laplacian}
L\phi(x)= \frac{m_U(x)}{m(x)}L_U \phi(x) + \frac{\phi(x)}{m(x)}(m(x)-m_U(x)),
\end{equation}
for any function $\phi$ satisfying $\phi=0$ on $ \textcolor{black}{{V_U}^c}$. The last term of the right hand side may be regarded as a boundary operator on 
$\partial U \textcolor{black}{=\{x \in V \mid N(x)\cap V_U \not=\emptyset \mbox{ and } N(x)\cap {V_U}^c} \\\textcolor{black}{\not=\emptyset\}. }$

\medskip

In what follows, we need the probability measure $\mu$ proportional to $m$ given by the normalization $\mu(A)=\frac{m(A)}{m(V)}$ for $A \subset V$.
We use the following estimate with the function $\mu(x)=\frac{m(x)}{m(V)}$ on $V$ later:
\begin{Lem}\label{estimate1} Let $U$ be a connected subgraph of $G$ and $f\in L^2(\mu)$. For any $\epsilon>0$ 
\[
\sum_{x\in V_U} f(x) \frac1{m(x)}\Delta_U f(x) \mu (x)\le \frac{\epsilon}{2} ||f||^2_{L^2(\mu)} + \frac{1}{\epsilon m(V)} {\mathcal E}(f,f).
\]
\end{Lem} 
{\it Proof.} By Fubini's theorem and {the} Cauchy-Schwarz inequality, we have
\begin{align*} 
& \sum_{x\in V_U} f(x) \frac1{m(x)}\Delta_U f(x)\mu(x) \\
& = m(V)\sum_{x\in V_U}\sum_{y\in V_U} 1_{\{\{x,y\}\in E_U\}} \frac1{m(x)m(y)}C_{x,y}(f(x)-f(y))f(x)\mu(x)\mu(y) \\
&\le m(V)\Big(\sum_{x\in V_U}\sum_{\{x,y\}\in E_U} f(x)^2 \frac1{m(x)m(y)}C_{x,y} \mu(x)\mu(y)\Big)^{1/2} \\
&\quad \times \Big(2\sum_{\{x,y\}\in E_U} (f(x)-f(y))^2 \frac1{m(x)m(y)}C_{x,y} \mu(x)\mu(y)\Big)^{1/2}.
\end{align*}
\medskip

Since $\sum_{\{x,y\}\in E_U} \frac1{m(x)m(y)}C_{x,y} \mu(y)\le \frac1{m(V)}\frac{1}{m(x)}
\sum_{\{x,y\}\in E_U} C_{x,y}\leq\frac1{m(V)}$, for any $\epsilon>0$, we see that the right-hand side is dominated by
$\epsilon^{1/2} \Big(\frac 2{\epsilon m(V)}\Big)^{1/2} ||f||_{L^2(\mu)} {\mathcal E}(f,f)^{1/2}$, which
does not exceed
$\frac{\epsilon}{2} ||f||_{L^2(\mu)}^2 +\frac {1}{\epsilon m(V)}  {\mathcal E}(f,f). $
\qed
\medskip

Here, we give a remark on the relationship between the weighted Laplacian $L$ and some stochastic processes: a reversible Markov chain 
$\{X_n\}_{n=0,1,2,\dots}$ and a Hunt process $\{Y_t\}_{t\in [0,\infty)}$. See \cite{Wo} for reversible Markov chains and see \cite{FOT} for Hunt processes and the theory of Dirichlet forms. 
\medskip

\begin{Remark} The Markov chain
$\{X_n\}$ can be determined by the transition matrix $P=(p(x,y))_{x,y\in V}$, where 
\begin{equation}\label{markov}
p(x,y)=C_{x,y}/m(x). 
\end{equation}
Thus for $A\subset V$, $P_x(X_1 \in A)=\sum_{x\in A} C_{x,y}/m(x)$. The Hunt process
$\{Y_t\}$ can be determined by transition semigroup $\{P_t\}$ associated with $L$ as its infinitesimal
generator. $\{P_t\}$ is a \textcolor{black}{\Erase{measure }}symmetric semigroup on $L^2(\mu)$; heuristically $\{P_t\}$ can be expressed by $P_t=e^{-tL}$ for any $t >0$.

\end{Remark}
\medskip

\subsection{Spectral gaps and Poincar\'e inequality}

 As in the second section, the graph metric between $x$ and $y$ is denoted by $d(x,y)$ and 
 $i(x)=\min \{\vert \Delta f(x)\vert \mid f: V_G\to \mathbb Z$ with 
 $\Delta f(x)\not=0\}$ is defined
 for every $x \in V_G$, due to the local finiteness of $G$. We fix a 
 base \textcolor{black}{\Erase{(reference)}}point $v_0$   
 and then take an exhaustive sequence $G_1\subset G_2 \subset \cdots$ of subgraphs of $G=(V_G,E_G)$ determined by 
$V_n=\{x \in V_G \mid d(v_0,x) \leq n\}$, $E_n=\{\{x,y\}\in E_G \mid x,y\in  V_n\}$ and $G_n=(V_n,E_n)$
for each positive integer $n$.

\medskip

We note that any divisor $D=\sum_{x \in V_{n-1}} \ell(x)i(x)1_{\{x\}}$ on the finite subgraph $G_n$ 
is viewed as the one on $G_j$ for any $j$ with $j>n$. In other words, $D=\sum_{x \in V_{n-1}} \ell(x)i(x)1_{\{x\}}$ can be consistently viewed as a divisor on $V_j$.

 \medskip

A formal linear combination 
$D=\sum_{x \in V_G} \ell(x)i(x)1_{\{x\}}$ of the family $\{1_{\{x\}}\mid x \in V_G\}$ of indicator functions
with 
at most countably many 
non-zero real coefficients $\{\ell(x)i(x)\mid x \in V_G\}$ is called a divisor on $G$.
Then, \Erase{due to Remark \ref{restrictionofdivisor}, $D\vert_{V_n}=\sum_{x \in V_{n}} \ell(x)i(x)1_{\{x\}}$ can not always be regarded 
as any divisor on $G_n$. On the other hand,}
the divisor $(D)_{n}=\sum_{x \in V_{n-1}} \ell(x)i(x) 1_{\{x\}}$ called 
{the} restriction of $D$ to $G_n$ is 
regarded as a divisor \textcolor{black}{vanishing outside $V_{n-1}$} on the finite graph $G_j$
for any integer $j$ with $j \geq n$. 
\Erase{We call 
$(D)_{n}$ restriction of $D$ to $V_n$.}
\medskip

Our objective is to  extend the Riemann-Roch theorem on finite graphs in \cite{BN} to the one on an infinite graph for 
divisor $D=\sum_{x \in V_G} \ell(x)i(x)1_{\{x\}}$ on $G$ satisfying $\sum_{x \in V_G} \vert \ell(x)\vert i(x)<\infty$.
To achieve this, it is necessary to discuss the convergence of $r_n(D)$ 
defined \Erase{as }in the second section on each $G_n$, 
at least in the case that ${\rm supp}[D]=\{x \in V_G \mid \ell(x)i(x)\not=0\}$ is a finite set. 
One of the key properties of the graphs for our aim is so-called {a} spectral  gap of $L$.   
We assume some conditions on the infinity of $G$ for taking our approach on the basis of spectral gap theory. 
\medskip

We denote the finite measure $m_{G_n}$ given in the previous subsection by $m_n$
and introduce the probability measure $\mu_{n}$ given by $\mu_n(A)=\dfrac{m_n(A)}{m_n(V_n)}$ for $A \subset V_n$.
For \textcolor{black}{\Erase{each }}$x \in 
 V_n$, $m_{n}(\{x\})$ is denoted briefly by $m_n(x)$. 

\medskip
In terms of a reversible Markov chain $\{X_n\}$ determined by the transition matrix defined in (\ref{markov}),
for any vertex $x$ with $d(x,v_0)=n$, we see that
$P_x(X_1\in V_n)=m_n(x)/m(x)$,
\[
P_x(X_1\in V_{n-1})=\frac1{m(x)} \sum_{y\in N(x)\cap V_{n-1}}C_{x,y}
\]
and these two identities imply $P_x(X_1\in V_n)\ge P_x(X_1\in V_{n-1})$. 
From those relationships among the probabilities, 
\[
\hbox{\Erase{$\rho_n\ge$}} P_x(X_1\in V_{n-1}^c)\ge \frac{m(x)-m_n(x)}{m(x)}.
\]
\noindent
%
\textcolor{black}{ \Erase{where
$S_n=\{x \in V_G \mid d(v_0,x)=n\}$
and}}
\textcolor{black}{ Define
\begin{equation}\label{rho_0}
\rho_n=\sup_{x\in 
S_n
} P_x(X_1\in V_{n-1}^c)
\end{equation}
for any positive integer $n$ where $S_n=\{x \in V_G \mid d(v_0,x)=n\}$. } \textcolor{black}{In the sequel, we assume that}

\begin{equation} \label{rho2} 
\limsup_{n\to\infty}\rho_n<1.
\end{equation} 
{B}y the assumption (\ref{rho2}), we can take $\rho<1$ such that $\limsup_{n\to\infty}\rho_n<\rho$.
Then we have 
\begin{equation}\label{rho3}
m(x)-m_n(x)\le \rho_n m(x)  \mbox{ and } m(x)\le \frac1{1-\rho} m_n(x)
\end{equation}

\noindent
for any $x\in 
S_n$ with sufficiently large $n$.
Since $m(x)=m_n(x)$ for any $x\in V_{n-1}$,
the first estimate and the finiteness of the measure $m$ imply
that
$m(V_n)-m_n(V_n)\to 0$
as $n\to \infty$. 

\medskip

\Erase{Here, we mention that asymptotic behavior of the reversible Markov chain $\{X_n\}$ is also controlled by the condition  (\ref{rho2}).
}
\medskip
In what follows, we denote $V_G$ by $V$, when\Erase{ the notation is preferably required for } it is preferable for the notation to require a shorter description.

\medskip
\begin{Lem}\label{poincare2} (\ref{rho2}) implies the inequality
\[
\sum_{x\in 
V_{n}, y\in N(x)\cap V_{n}^c} f(x)^2 C_{x,y} \le \rho_n m(V) ||f||_{L^2(\mu)}^2.
\]
\noindent
on real-valued function $f$ on $V_{n}^c$ for sufficiently large $n$.
\end{Lem}
{\it Proof.} \quad The left hand side equals $\sum_{x\in 
S_n}
f(x)^2(m(x)-m_n(x))$ and
 (\ref{rho3}) shows that this value is dominated by
$\rho_n \sum_{x\in S_{n}}f(x)^2 m(x)$
for sufficiently large $n$.
 \qed
\medskip

The following assertion is crucial to establish Riemann-Roch theorem 
in our approach.

\medskip

\begin{Thm}\label{poincare} There exists a positive constant $A$ with $A <1$ such that 
\[
\limsup_{n\to\infty}\rho_n<A
\]
implies the existence of a positive constant $C$ such that
\begin{equation}\label{poincare ineq}
||f||^2_{L^2(\mu)}\le C{\mathcal E}(f,f),
\end{equation}
for any $f\in L^2(\mu)$ with $(f,1)_{L^2(\mu)}=0$. 
\end{Thm} 


\medskip

\noindent

\medskip

{\it Proof.}  
We first divide the sum of the left-hand side of (\ref{poincare ineq}) into two terms as
\begin{equation}\label{division}
||f||^2_{L^2(\mu)}=\sum_{x\in V_n} f(x)^2\mu(x) + \sum_{x\in V_n^c} f(x)^2\mu(x).
\end{equation}

\noindent
The first term can be written as
\[
\sum_{x\in V_n} f(x)^2\mu(x) =\frac{m_n(V_n)}{m(V)}||f||^2_{L^2(\mu_n)}+\frac1{m(V)}\sum_{x\in V_n} f(x)^2(m(x)-m_n(x)).
\]

\noindent
Since $G_n$ is a finite graph, there exists a positive constant $\lambda_n$ such that 
\begin{align*}
\sum_{x\in V_n} f(x)^2\mu(x) &\le \frac{m_n(V_n)}{m(V)\lambda_n}{\mathcal E}^{(n)}(f,f)+ \frac{m_n(V_n)}{m(V)}(f,1)_{L^2(\mu_n)}^2 +\nu_n(f)\\
&\le \frac{m_n(V_n)}{m(V) \lambda_n}{\mathcal E}(f,f)+\frac{m_n(V_n)}{m(V)} (f,1)_{L^2(\mu_n)}^2 +\nu_n(f),
\end{align*}
where ${\mathcal E}^{(n)}(f,f)={\mathcal E}_{G_n}(f,f)$
and $\nu_n(f)=\frac1{m(V)}\sum_{x\in V_n} f(x)^2(m(x)-m_n(x))$. 
Since $m(x)=m_n(x)$ for any $x\in V_{n-1}$, by Lemma \ref{poincare2} we obtain
\[
\nu_n(f)\le \rho_n ||f||^2_{L^2(\mu)}
\]
and by $(f,1)_{L^2(\mu)}=0$, 
\[
(f,1)_{L^2(\mu_n)}^2=\Big(\frac{m(V)}{m(V_n)}\sum_{x\in V_n^c}|f(x)|\mu(x)\Big)^2
\le (\frac{m(V)}{m(V_n)})^2 \mu(V_n^c)  ||f||^2_{L^2(\mu)}.
\]
Combining these estimates, we have
\begin{equation}\label{inner}
\sum_{x\in V_n} f(x)^2\mu(x) \le \frac{m_n(V_n)}{m(V)}\lambda_n^{-1}{\mathcal E}(f,f)
+\alpha_n ||f||^2_{L^2(\mu)},
\end{equation}
where $\alpha_n= \rho_n + (\frac{m(V)}{m(V_n)})^2 \mu(V_n^c) $. 
\medskip

Next{,} we estimate the second term of the right hand side of (\ref{division}) by
taking the function $e^{a r(x)}$ on $V$ \Erase{defined by }with $r(x)=d(x,v_0)$ and {an} arbitrarily fixed $a>0$.
\Erase{For the function $e^{a r(x)}$, we have} {Since} 
\begin{align*}
\Delta e^{a r(x)} & = \sum_{y\in N(x), r(y)=r(x)+1} C_{x,y} (e^{a r(x)}-e^{a r(y)}) +\sum_{y\in N(x), r(y)=r(x)-1} C_{x,y} (e^{a r(x)}-e^{a r(y)}) \\
&= -(e^a-1)e^{a r(x)}m_+(x) + (1-e^{-a})e^{a r(x)}m_-(x),
\end{align*}
where 
\[
m_+(x)=\sum_{y\in N(x), r(y)=r(x)+1} C_{x,y} \quad \mbox{ and }\quad m_-(x)=\sum_{y\in N(x), r(y)=r(x)-1} C_{x,y},
\]

\noindent
\Erase{In other words,} we have
$L e^{a r(x)}=  -(e^a-1)e^{a r(x)}\frac{m_+(x)}{m(x)} + (1-e^{-a})e^{a r(x)}\frac{m_-(x)}{m(x)}.$
\Erase{We note}{Noting} that  
\[
\frac{m_-(x)}{m(x)}=P_x(X_1\in V_{n-1})
\]
for any $x\in S_n$ with $n\geq 1$, \Erase{and } {we have} another expression on $\rho_n$, \Erase{given as}
$\rho_n=\sup\{1-\frac{m_-(x)}{m(x)}\mid x\in S_n\}$, 
\Erase{ shows that }
{and thus}
\[
\frac{m_+(x)}{m(x)} \le 1-\frac{m_-(x)}{m(x)}\le \rho_n
\]

\noindent
for any $x\in S_n$ with $n\geq 1$.

From these observations, we can derive that
\begin{align*}
e^{-a r(x)}L e^{a r(x)} &=  -(e^a-1)\frac{m_+(x)}{m(x)} + (1-e^{-a})\frac{m_-(x)}{m(x)} \\
&\ge -(e^a-1)(1-\frac{m_-(x)}{m(x)}) - (1-e^{-a})(1-\frac{m_-(x)}{m(x)})+1-e^{-a} \\
&\ge -(e^a-e^{-a})\rho_n +1-e^{-a}.
\end{align*}
For the sequence $\{\eta_n\}$
defined by $\eta_n= -(e^a-e^{-a})\rho_n +1-e^{-a}\quad (n=1,2,\dots)$
, we
may assume that $\eta_n>0$ for any $n\ge n_0$
with some positive integer $n_0$, more specifically,

\begin{equation}\label{rho condition1}
\rho_n<\frac{1-e^{-a}}{e^a-e^{-a}} \quad (n\ge n_0). 
\end{equation}
\medskip

Let us apply (\ref{laplacian}) to the subgraph 
$U=G_n^c$ 
determined by taking $V_n^c$ as $V_U$ 
and the function $\phi$ defined by

$$\phi(x)=\begin{cases} e^{a r(x)}\; &(x\in V_{n-1}^c),\\
0 \; &(x\in V_{n-1}).
\end{cases}$$
Then, we have
$
e^{-a r(x)} L\phi(x)=\frac{m_{G_n^c}(x)}{m(x)}e^{-a r(x)} L_{G_n^c}\phi(x)+ \frac1{m(x)}\sum_{y\in N(x)\cap V_{n-1}}C_{x,y},
$
which implies
\begin{align}\label{bdry}
\sum_{x\in V_n^c} f(x)^2\mu(x) &=\frac1{\eta_n}\sum_{x\in V_{n-1}^c} f(x)^2e^{-a r(x)} L\phi(x)\mu(x) - \sum_{x\in S_n}f(x)^2\mu(x) \notag \\
 &=\frac 1{\eta_n m(V)}\sum_{x\in V_{n-1}^c} f(x)^2 e^{-a r(x)} L_{G_{n-1}^c}\phi(x) m_{G_{n-1}^c}(x) \notag \\
& \qquad +\frac1{m(V)}\Big\{ \frac1{\eta_n}\sum_{x\in V_{n-1}^c} f(x)^2\sum_{y\in N(x)\cap V_{n-1}}C_{x,y}-\sum_{x\in S_n}f(x)^2 m(x)\Big\}.
\end{align}
We can regard the second term in the right-hand side as the one involving a boundary operator.
We estimate the first term in the right hand side of (\ref{bdry}). 
By the symmetry of $L_{G_{n-1}^c}$  with respect to the 
inner product of $L^2(m_{G_{n-1}^c})$, 
\begin{align*}
\sum_{x\in V_{n-1}^c} f(x)^2 e^{-a r(x)}  L_{G_{n-1}^c}\phi(x) m_{G_{n-1}^c}(x) 
&=\sum_{x\in V_{n-1}^c}  e^{ar(x)} L_{G_{n-1}^c}(f(x)^2 e^{-ar(x)}) m_{G_{n-1}^c}(x) \\ 
&=\sum_{x\in V_{n-1}^c}  e^{ar(x)} \Delta_{G_{n-1}^c}(f(x)^2 e^{-ar(x)}). 
\end{align*}
By direct calculations, we have 
\[
 \Delta_{G_{n-1}^c}(f(x)^2 e^{-a r(x)}) =  e^{-a r(x)} \Delta_{G_{n-1}^c}(f(x)^2)+ f(x)^2\Delta_{G_{n-1}^c}(e^{-a r(x)}) 
-\Gamma(f(x)^2, e^{-a r(x)}),
\]
where $\Gamma(u,v)(x)=\sum_{y\in N(x)\cap V_{n-1}^c}C_{x,y}(u(x)-u(y))(v(x)-v(y))$. 
It is easy to see that
\begin{align}
\Delta_{G_{n-1}^c}(f(x)^2) &=2f(x) \Delta_{G_{n-1}^c}f(x) -\Gamma(f,f)
\intertext{and}
\Delta_{G_{n-1}^c}(e^{-a r(x)}) &= e^{-a r(x)}(1-e^{-a})m_+^{G_{n-1}^c}(x)
-(e^{a}-1)e^{-a r(x)}m_-^{G_{n-1}^c}(x)\nonumber \\ 
&\le e^{-a r(x)}(1-e^{-a})\rho_n m(x),
\end{align}
where $m_+^{G_{n-1}^c}(x)=\sum_{r(y)=r(x)+1, y\in V_{n-1}^c}
C_{x,y}$. 
We also see that
\[
|\Gamma(f(x)^2, e^{-a r(x)})|\le (e^a-1)e^{-a r(x)}\sum_{y\in N(x)\cap V_{n-1}^c}C_{x,y}|f(x)+f(y)||f(x)-f(y)|.
\]
From these observations, we have
\begin{align*}
& \frac 1{\eta_n m(V)}\sum_{x\in V_{n-1}^c}  e^{a r(x)} \Delta_{G_{n-1}^c}(f(x)^2 e^{-a r(x)}) \\
&\le \frac 1{\eta_n m(V)}\Big\{\sum_{x\in V_{n-1}^c} (2f(x) \Delta_{G_{n-1}^c}f(x) -\Gamma(f,f)) \\ 
& \qquad + (1-e^{-a})\rho_n\sum_{x\in V_{n-1}^c} f(x)^2 m(x) \\ 
&\qquad+(e^a-1)\sum_{x\in V_{n-1}^c}\sum_{y\in N(x)\cap G_{n-1}^c}C_{x,y}|f(x)+f(y)||f(x)-f(y)|\Big\} \\ 
&\le \frac 1{\eta_n}\Big\{\sum_{x\in V_{n-1}^c} 2f(x) \frac1{m(x)}\Delta_{G_{n-1}^c}f(x)\mu(x) \\ 
&\qquad +(1-e^{-a})\rho_n\sum_{x\in V_{n-1}^c} f(x)^2 \mu(x) \\ 
&\qquad+ (e^a-1)\sum_{x\in V_{n-1}^c}\sum_{y\in N(x)\cap V_{n-1}^c}C_{x,y}|f(x)+f(y)||f(x)-f(y)|\Big\}.
\end{align*} 
Lemma \ref{estimate1} implies that the first term in the right hand side is dominated by
\[
\frac {\epsilon}{\eta_n} ||f||^2_{L^2(\mu)}+  \frac{2 }{\epsilon m(V)\eta_n}{\mathcal E}(f,f)
\]
for any $\epsilon>0$. The second term is dominated by 

\[
\frac {(1-e^{-a})\rho_n}{\eta_n} ||f||^2_{L^2(\mu)}.
\]

As for the third term, by {the} Cauchy-Schwarz inequality{,} we have
\begin{align*}
&\sum_{x\in V_{n-1}^c}\sum_{y\in N(x)\cap V_{n-1}^c}C_{x,y}|f(x)+f(y)||f(x)-f(y)|\\
&\quad= \sum_{x\in V_{n-1}^c}\sum_{y\in V_{n-1}^c} 1_{\{x,y\}\in E_{V_{n-1}^c}} C_{x,y}|f(x)+f(y)||f(x)-f(y)|
\\
&\quad\le \big(\sum_{x\in V_{n-1}^c}\sum_{y\in V_{n-1}^c} 1_{\{x,y\}\in E_{V_{n-1}^c}} C_{x,y}|f(x)+f(y)|^2\big)^{1/2}
(2\mathcal{E}(f,f))^{1/2}.
\end{align*}
Fubini's theorem implies 
\[
\sum_{x\in V_{n-1}^c}\sum_{y\in V_{n-1}^c} 1_{\{x,y\}\in E_{V_{n-1}^c}} C_{x,y}|f(x)+f(y)|^2 \le 
4 m(V)||f||^2_{L^2(\mu)}.
\]
Thus, the third term does not exceed 
\[
\frac {(e^a-1)\epsilon}{\eta_n}||f||^2_{L^2(\mu)}+ \frac {4(e^a-1)m(V)}{\epsilon \eta_n}{\mathcal{E}}(f,f).
\]

%
%
%
We next estimate the term with the boundary operator  in (\ref{bdry}) as
\begin{align*}
&\frac1{m(V)}\Big\{ \frac1{\eta_n}\sum_{x\in V_{n-1}^c} f(x)^2\sum_{y\in N(x)\cap V_{n-1}}C_{x,y}-\sum_{x\in S_n}f(x)^2 m(x)\Big\}\\
&\quad=\frac1{m(V)}\sum_{x\in S_n} f(x)^2\big(\frac1{\eta_n}\sum_{y\in N(x)\cap V_{n-1}}C_{x,y}- m(x)\big) \\
&\quad\le \frac1{m(V)}\sum_{x\in S_n} f(x)^2( \frac1{\eta_n}-1)m_-(x) \\
&\quad\le( \frac1{\eta_n}-1) ||f||^2_{L^2(\mu)}.
\end{align*}
Combining these estimates, we have
\begin{equation}\label{outer}
\sum_{x\in V_n^c} f(x)^2\mu(x)\le   \beta(n,\epsilon) ||f||^2_{L^2(\mu)} + \gamma(\epsilon,n){\mathcal E}(f,f),
\end{equation}
where 
$
\beta(\epsilon,n)= \frac {(e^a+1)\epsilon+2(1-e^{-a})\rho_n}{2\eta_n}+ \frac1{\eta_n}-1
$
and 
$
\gamma(\epsilon,n)=\frac{2 m(V)^{-1}+4(e^a-1)m(V)}{\epsilon \eta_n}.
$
From (\ref{inner}) and (\ref{outer}), we can derive

\[
||f||^2_{L^2(\mu)}\le \big(1+\lambda_n^{-1}+\gamma(\epsilon,n)\big) {\mathcal E}(f,f)+\big(\alpha_n+\beta(\epsilon,n)\big)||f||^2_{L^2(\mu)}.
\]

\medskip

%
 
 \medskip
{Here, we} consider a condition on $\rho_n$ ensuring that $\alpha_n+\beta(\epsilon,n)<1$. 
Let $B(a),B'(a)$ be the smaller and larger solutions of the quadratic equation
\[
(e^a-e^{-a})t^2-2(e^a-2e^{-a}+1)t+1-2e^{-a}=0,
\]
respectively. Then $\rho_n<B(a)$ implies
\[
\rho_n + \frac{(1-e^{-a})\rho_n}{\eta_n}+ \frac1{\eta_n}-1<1.
\]

We can take the maximum of $B(a)$ subject to $a>\log 2$ as the positive constant $A$ in the \textcolor{black}{statement of Theorem \ref{poincare}}. 
In fact,  $B(a)>0$ for $a>\log 2$ and under the condition
 that $\limsup_{n\to\infty}\rho_n<\max \{B(a)\mid a>\log 2\}$,
we can easily take a positive real $a$ such that the condition (\ref{rho condition1}) is satisfied for sufficiently large $n$. (This is because $B(a)<\frac{1-e^{-a}}{e^a-e^{-a}}<B'(a)$ for any $a$). Consequently, there exist $n_G$ and $\epsilon_0>0$ such that
$\alpha_{n_G}+\beta(\epsilon_0, n_G)<1$. Therefore, we have  
\begin{equation}\label{n_G}
||f||^2_{L^2(\mu)}\le \frac{1+\lambda_{n_G}^{-1}+\gamma(\epsilon_0, n_G)}{1-\alpha_{n_G}-\beta(\epsilon_0,n_G)} {\mathcal E}(f,f). 
\end{equation}
%
%
\qed

\begin{Remark}\label{Remark 3.6'} \textcolor{black}{When we apply  this theorem for a function $f$ on subgraph $G_n=(V_n,E_n)$,
i.e., a function with its support contained in $V_{n}$ and $(1,f)_{L^2(\mu_{n})}=0$, 
we need to restrict our attention to such sufficiently large $n$ that $\rho_n < A$ holds. This point will be used when we apply Lemma \ref{L2-bdd} below in Section 4.   
}
\end{Remark}

\begin{Remark} By a numerical calculation, one sees that $\textcolor{black}{A=}\max \{B(a)\mid a>\log 2\}\approx 0.0569$.
\end{Remark}
%

For a connected subgraph $U$ of $G$, we introduce
\[
\lambda_U=\inf\{\tilde{\mathcal E}_{U}(f,f)\; \mid\; f\in L^2(\mu_U), (f,1)_{L^2(\mu_U)}=0\text{ and }\Vert f\Vert_{L^2(\mu_U)}=1\}
\] 
called a spectral gap of $L_U$, where $\tilde{\mathcal{E}}_U(u,v)=m_U(V_U)^{-1}\mathcal{E}_U(u,v)$. 
Similarly to the relationship between the Laplace operator $L_U$ on $L^2(m_U)$ and $\mathcal E_U$,
we note that $(L_U u,v)_{L^2(\mu_U)}=\tilde{\mathcal E}_U(u,v)$. It is well known that if $U$ is a finite graph, then $\lambda_U>0$
(cf. \cite{FN}). 
The spectral gap $\lambda_n$ of $L_{G_n}$ is 
taken in the first paragraph of the proof of Theorem \ref{poincare},
the assertion of which shows the positivity of the spectral gap of $L$. 
\begin{Cor} Under the assumption of Theorem \ref{poincare}, we have
\[
\lambda_G\ge \frac{1}{Cm(
V_G)
}.
\]
\end{Cor}
%
\begin{Lem}\label{spectral gap}
Let 
$\lambda_n$ be the spectral gap of $L_{G_n}$
and assume $\lambda_{G}>0$. 
If $\limsup_{n\to\infty}$\linebreak
$\frac{\rho_n}{1-\rho_n}<\lambda_G$, then
$\liminf_{n\to\infty}\lambda_n>0$. Moreover, if $\lim_{n\to\infty} \rho_n=0$, then 
\[
\liminf_{n\to\infty} \lambda_n\ge \lambda_G.
\]
\end{Lem}
{\it Proof.} Since $G_n$ is a finite graph, $\lambda_n$ is the smallest 
non-zero
eigenvalue of $L_{G_n}$
and any eigenfunction $\psi_n$ associated with $\lambda_n$ satisfies $(1,\psi_n)_{L^2(\mu_n)}=0$.
We assume that the eigenfunction is normalized as
$||\psi_n||_{L^2(\mu_{n})}=1$ beforehand.
We extend this function to the one defined on $V$
taking identically zero on $V_n^c$ and denote it again by $\psi_n$. 
Let $\overline{\psi_n}=\psi_n-(1,\psi_n)_{L^2(\mu)}$ and we can first show that $||\overline{\psi_n}||_{L^2(\mu)}\ge m(V_n)m(V)^{-1}$. 
In fact, the identities
\begin{align*}
||\overline{\psi_n}||_{L^2(\mu)} &= \sum \big(\psi_n(x)-(1,\psi_n)_{L^2(\mu)}\big)^2\mu(x) \\
&= \sum \psi_n(x)^2\mu(x)-(1,\psi_n)_{L^2(\mu)}^2 \\
&= m(V)^{-1}\sum \psi_n(x)^2(m(x)-m_n(x))+m_n(V_n)m(V)^{-1}\\
&\quad -\Big(m(V)^{-1}\sum \psi_n(x)(m(x)-m_n(x))\Big)^2 \\
\end{align*}
and {the} Cauchy-Schwarz inequality yield a lower bound 
\begin{align*}
& m(V)^{-1}\sum \psi_n(x)^2(m(x)-m_n(x))+m_n(V_n)m(V)^{-1}\\
&\quad-m(V)^{-2}(m(V)-m_n(V_n))\sum \psi_n(x)^2(m(x)-m_n(x)) \\
&\ge m_n(V_n)m(V)^{-1}
\end{align*}

\noindent
of the right-hand side.
Since $m(x)-m_n(x)\le \frac{\rho_n}{1-\rho_n} m_n(x)$, 
\begin{align*}\label{lower}
	{\mathcal E}(\overline{\psi_n},\overline{\psi_n})&= {\mathcal E}_{G_n}(\overline{\psi_n},\overline{\psi_n})+ \sum_{x\in G_n, y\in V_n^c} C_{x,y}\psi_n(x)^2 \\
&=m_n(V_n)\lambda_n+\sum_{x\in V_n}\psi_n(x)^2(m(x)-m_n(x)) \\
&\le m_n(V_n)\lambda_n+\frac{\rho_n}{1-\rho_n}m_n(V_n).
\end{align*}
As a result, it turns out that
\[
m_n(V_n)^{-1}m(V)\lambda_G\le \lambda_n+\frac{\rho_n}{1-\rho_n}.
\]
\qed
%
\begin{Lem}\label{L2-bdd} 
If  $\limsup_{n\to\infty} \frac{\rho_n}{1-\rho_n}<\lambda_G$ as assumed in Lemma \ref{spectral gap}, then there exists some positive constant $K$ such that any sequence $\{g_n\}$  of functions $g_n$ on $V_n$ satisfying
$(1,g_{n})_{L^2(\mu_{n})}=0$ for sufficiently large $n$ yields
\[
||R^{(n)}g_n||_{L^2(\mu_{n})} \le K||g_n||_{L^2(\mu_{n})}
\]
for sufficiently large $n$, where $R^{(n)}$ stands for the $0$-order resolvent of $L_{G_n}$ on $L^2(\mu_n)$.
\end{Lem}
{\it Proof.} For the transition operator $P^{(n)}_t$ associated with the generator $L_{G_n}$ on $L^2(\mu_n)$, the positivity of spectral gaps implies (cf.\cite{Wa})
\[
||P^{(n)}_t g_{n}||_{L^2(\mu_{n})}\le e^{-\lambda_n t} ||g_{n}||_{L^2(\mu_{n})}.
\]
Note that $R^{(n)}g_{n}(x)=\int_0^{\infty}P^{(n)}_t g_{n}(x)dt$ and 
\[
R^{(n)}g_{n}(x)=\sum_{k=1}^{\infty} R^k_1g_{n}(x), 
\]
where $R_1g_{n}(x)=\int_0^{\infty}e^{-t}P^{(n)}_t g_{n}(x)dt$ and $R^k_1$ denotes the operator given as the $k$-times iteration of $R_1$ 
(cf. \cite{FOT}). 
Lemma \ref{spectral gap} implies that there exists  $\delta>0$ such that  
$\lambda_n\ge \delta$ for all $n$.   
By applying Jensen's inequality and Fubini's theorem, we obtain
\[
||R_1g_{n}||_{L^2(\mu_{n})}^2\le \frac 1{2\lambda_n+1}||g_{n}||_{L^2(\mu_{n})}^2\le \frac1{2\delta+1}||g_{n}||_{L^2(\mu_{n})}^2
\]
for sufficiently large $n$. 
The desired result follows from this estimate. 
\qed

\section{Proof of Riemann-Roch theorem on an infinite graph.}

In this section, we establish a Riemann-Roch theorem on the connected infinite graph $G=(V_G,E_G)$
with local finiteness and finite volume as in the last section,
by applying $L^2$-boundedness of the $0$-order resolvent derived from spectral gap theory to
a sequence of functions in the images of the Laplace operator.
For that purpose, we first take a divisor $D=\sum_{x\in V_{n-1}}\ell(x)i(x) 1_{\{x\}}$ on some subgraph $G_n=(V_n,E_n)$ 
of $G=(V_G,E_G)$ as given in the last section
and denote $\Delta_{G_n} $ and $L_{G_n} $ by $\Delta_n$ and by $L_n$, respectively. 
The 
equivalence
 between divisors $D'=\sum_{x\in V_{n}}\ell'(x)i(x) 1_{\{x\}}$ and $D''=
\sum_{x\in V_{n}}\ell''(x)i(x) 1_{\{x\}}$  
with $\mathbb Z$-valued functions $\ell'$ and $\ell''$ is defined  by
$D''=D' + \Delta_n f$ on $V_{n}$ for some $\mathbb Z$-valued function \textcolor{black}{$f$} with ${\rm supp}[f]\subset V_{n}$.
This relationship
will be
denoted by $D' \overset {n}{ \sim} D''$ and 
will be called $n$-equivalence.  
The family
of total orders on $V_{n}$ is denoted by $\mathcal O_n$.
 \medskip

For  a divisor $D$ on {the} finite graph $G_n=(V_n,E_n)$, 
$r_n(D)$ is defined on the finite subgraph $G_{n}$, 
by replacing ``$\sim$" in the second section with the $n$-equivalence ``$\overset {n}{ \sim}$", 
more specifically given by
$$r_n(D)
=\Big( \min_{D'\overset {n}{ \sim} D, O_n \in \mathcal O_n}\mbox{\rm deg}^{+}(D' - \nu_{O_n})\Big)-i_{(G_n,C_n)}.$$

When a divisor $D=\sum_{x \in V_G} \ell(x)i(x)1_{\{x\}}$ satisfying $\sum_{x \in V_G} \vert \ell(x)\vert i(x)<\infty$ is given,
\textcolor{black}{we denote $\sum_{x \in V_G,\ell(x)>0} \ell(x)i(x)$ by ${\rm deg}^+(D)$
and $-\sum_{x \in V_G,\ell(x)<0}\ell(x)i(x)$ by ${\rm deg}^-(D)$.} 
The family $\{(D)_{n}\}$ of divisors \Erase{has a consistency}\textcolor{black}{is consistent} in the sense that
$((D)_j)_n=(D)_{n}$ whenever $j>n$.
Later, by 
taking control over
the sequence $\{O_j\}$ of total orders each of which is taken in the minimization of $r_j((D)_{n})
=\min_{D'\overset {j}{ \sim} (D)_n, O_{j} \in \mathcal O_{j}}{\rm deg}^+(D'-\nu_{O_j})-i_{(G_j,C_j)}$ 
on every subgraph $G_j$ 
satisfying $V_j \supset V_{n-1}\supset {\rm supp}[(D)_{n}]$, we 
facilitate successive procedures of taking limits as $j\to \infty$ and $n \to \infty$ in $r_j((D)_{n})$.
\Erase{so that those limits are taken along an identical subsequence
of positive integers.}
\textcolor{black}{In fact, thanks to Proposition \ref{fundamental}, we see that
$\vert r_j((D)_{n})- r_j((D)_{n'})\vert$ is sufficiently small as long as $n,n'$ are both large enough,
independently of $j$ with $j > \max\{n,n'\}$. This is because ${\rm deg}^+((D)_{n}-(D)_{n'})+{\rm deg}^-((D)_{n}-(D)_{n'})$
is sufficiently small as long as $n,n'$ are both large enough.}
 \medskip

\Erase{
In the following two lemmas, we focus on a divisor $D$ satisfying ${\rm supp}[D] \subset V_{n-1}$ for some $n$
and a minimizer for $r_n(D)$ is denoted by $f_n$. }
%
\textcolor{black}{
In the following two lemmas, we focus only on a divisor $D$ satisfying
${\rm supp}[D] \subset V_{j_0}$ for some $j_0$
and the case that \Erase{${\rm supp}[D] \subset V_{n-1}$ and }$\rho_n< \min\{A,\frac{\lambda_G}{\lambda_G+1}\}$ is 
satisfied for sufficiently large $n>j_0$.
We note that, in this case, the uniform boundedness of the family
$\{R^{(n)}\}$ of operators is derived from the proof of Lemma \ref{L2-bdd}.
We denote
a minimizer for $r_n(D)$ by $f_n$.}
\textcolor{black}{
To be more precise,
$f_n$ is an integer-valued function on $V_{n}$ such that $D'= D + \Delta_{n}f_n$ attains the minimum }
 
 \begin{equation*}
r_{n}(D)
=\Big( \min_{D'\overset {n}{ \sim} D, 
O_{n} \in \mathcal O_{n}}
{\rm deg}^{+}(D' - \nu_{O_{n}})\Big)-i_{(G_{n},C_{n})},
\end{equation*}

 \medskip
 \noindent
 with some  total order $O_{n} \in \mathcal O_{n}$.
 
  \medskip
\textcolor{black}{ For such a divisor $D$ 
the divisor $(D)_j$ is equal to $D$ for any $j$ with $j\geq n$.}
\textcolor{black}{To show the convergence of a subsequence of $\{ r_{j}(D) \}$, we will need an upper bound of $\vert f_n\vert$ on $V_n$.
In fact, the boundedness enables us to take the real-valued harmonic extension $h_j^{(n)}$ of $f_n$ on $V_j \setminus V_{n}$ as seen after the
following lemma, which yields a sequence of integer-valued functions used for deriving the convergence of a subsequence later.}
 \medskip

 \begin{Lem}\label{subsequence}
\Erase{For any positive integer $n$, } For sufficiently large integer $n$, $r_{n}(D)$ 
admits a minimizer $f_{n}$  satitsfiying
$\Vert f_{n}\Vert_{L^2(\mu_{n})}^2 \leq 4 m(V)K \big(({\rm deg^+}(D)+{\rm deg^-}(D)+m_n(V_n))^2/\min\{m(x)\mid x \in V_{n}\}\big)$,
where $K$ stands for the constant given in Lemma \ref{L2-bdd}. In particular,
$$\max_{x \in V_{n}}\vert f_{n}(x)\vert \leq 2\sqrt{m(V)K}({\rm deg^+}(D)+{\rm deg^-}(D)+1)/\min\{m(x)\mid x \in V_{n}\}$$

\noindent
for sufficiently large $n$.
 \end{Lem}
 
 {\it Proof.} 
Any $\mathbb Z$-valued function $f^{(n)}_{D}$ on $V_n$ minimizing ${\rm deg}^+(D + \Delta_{n}f^{(n)}_{D})$ 
satisfies ${\rm deg}^+(D + \Delta_{n}f^{(n)}_{D})\leq {\rm deg}^+D$. Accordingly, Lemma \ref{deg} (ii) shows that 
${\rm deg}^-(D + \Delta_{n}f^{(n)}_{D})\leq {\rm deg}^-D$, which implies $\vert \Delta_{n}f^{(n)}_{D}(x) \vert$ 
does not exceed $2({\rm deg^+}(D)+{\rm deg^-}(D))$ for any $x\in V_n$ and $f^{(n)}_{D}$
admits the estimate
$\Vert \frac{1}{m}\Delta_{n}f^{(n)}_{D}\Vert_{L^2(\mu_{n})}^2 \leq 4 m(V)({\rm deg^+}(D)+{\rm deg^-}(D))^2/\min\{m(x)\mid x \in V_{n}\}$. Since $\vert\nu_{O_{n}}(x)\vert
 \leq m(x)$ for every $x \in V_{n}$ and $O_n \in \mathcal O_n$,
any minimizer $f_{n}$ for $r_n(D)={\rm deg}^+(D + \Delta_{n}f_{n}- \nu_{O_{n}})$, which is viewed as  $f^{(n)}_{D+\nu_{O_{n}}}$, satisfies
\begin{align*}
\Vert \frac{1}{m}\Delta_{n}f_{n}\Vert_{L^2(\mu_{n})}^2 
& \leq 4m(V)({\rm deg^+}(D)+{\rm deg^-}(D)+m_n(V_n))^2/\min\{m(x)\mid x \in V_{n}\}.
\end{align*}

\noindent
Thanks to Lemma \ref{L2-bdd},
these inequalities imply 
$\Vert R^{(n)} L_{n} f_{n}\Vert_{L^2(\mu_{n})}^2
\leq 4 m(V)\linebreak K\big(({\rm deg^+}(D)
+{\rm deg^-}(D)+m_n(V_n))^2/\min\{m(x)\mid x \in V_{n}\}\big)$ .
\medskip

The uniqueness of the solution to the
Poisson equation 
on connected finite graph{s} 
shows that $y,z\in V_n$ implies $R^{(n)}L_{n}f_{n}(y)-R^{(n)}L_{n} f_{n}(z) \in \mathbb Z$. This and the uniqueness up to difference given by constants 
allow us to take such a real constant $c$ with
$0\leq c < 1$ that $R^{(n)}L_{n} f_{n}+c$ is not only $\mathbb Z$-valued but a minimizer for $r_{n}(D)$.
We focus only on this minimizer obtained by this procedure and denote it by $f_n$, then $f_n$ admits 
the estimate in the assertion and  $\max_{x \in V_{n}}\vert f_{n}(x)\vert \leq 2\sqrt{ m(V)K}\big(({\rm deg^+}(D)+{\rm deg^-}(D)+1)/\min\{m(x)\mid x \in V_{n}\}\big)$  is satisfied 
for sufficiently large $n$. In fact,  this follows from $\vert f_{n}(x)\vert^2m(x) \leq 4 m(V)K\big(({\rm deg^+}(D)+{\rm deg^-}(D)+1)^2/\min\{m(x)\mid x \in V_{n}\}\big)$ for any $x \in V_{n}$
as long as $n$ is sufficiently large.\qed

  \medskip
\begin{Lem}{\label{(I)}} \Erase{For each $n$ }For sufficiently large integer $n$, the minimizer $f_{n}$  in the last lemma for $r_{n}(D)$ admits a sequence $\{f^{(n)}_{j}\}_{j>n}$
 of functions such that

\begin{description}
\item[(i)] $ f_{n}= f^{(n)}_{j}$ on $V_{n}$, 
\item[(ii)] $\max_{x \in V_{j}}\vert f_{j}^{(n)}(x)\vert 
\leq 2\sqrt{m(V)K}({\rm deg^+}(D)+{\rm deg^-}(D)+1)/\min\{m(x)\mid x \in V_{n}\} \\
\qquad\qquad\qquad\qquad
+1$,
\item[(iii)] $\lim_{n\to \infty}\sup_{j > n}\Vert m_{j}(\cdot)^{-1}(\Delta_{j} f^{(n)}_{j}-D) \Vert_{L{\phantom{}^1(V_{j}\setminus V_{n}};\mu)}=0$.
\end{description}
\end{Lem}
  \medskip
  
  {\it Proof.} We take the harmonic extension $h^{(n)}_{j}$ \textcolor{black}{of $f_n$} on $V_{j}\setminus \textcolor{black}{V_{n}}$  \textcolor{black}{determined by $\Delta_j h^{(n)}_{j} =0$ on $V_j \setminus V_n$ and $h^{(n)}_{j}=f_{n}$ on $V_{n}$.}
The function $g^{(n)}_{j}$ defined as the integer part of $h^{(n)}_{j}$ admits
the estimate 

$$\frac{\vert \Delta_{j} g^{(n)}_{j}(x)\vert}{m(x)}\big\vert_{V_{j}\setminus V_{n}}\leq 
21_{V_{j}\setminus V_{n}}$$

\noindent
on $\Delta_{j} g^{(n)}_{j}(x)=  \sum_{y \in N(x)\cap V_j}C_{x,y}(g^{(n)}_{j}(x)-g^{(n)}_{j}(y))$.
\medskip

In fact, since
$\max_{x \in V_{j}\setminus V_{n-1}} \vert  h^{(n)}_{j}(x) - g^{(n)}_{j}(x)\vert \leq 1$, 
we see 
\medskip

$$ C_{x,y}\vert g^{(n)}_{j}(x)-h^{(n)}_{j}(x) - g^{(n)}_{j}(y)-h^{(n)}_{j}(y)\vert\leq 2C_{x,y}$$

\noindent
for any $x\in V_j\setminus V_n$ and $y\in N(x)$.
By combining this with the harmonicity of $h^{(n)}_{j}$ on $V_{j}\setminus V_{n-1}$, 
equivalently $x \in V_{j}\setminus V_{n}\Rightarrow \sum_{y \in N(x)\cap V_j}C_{x,y}
(h^{(n)}_{j}(x)-h^{(n)}_{j}(y))=0$, we see that 
$\frac{\vert\sum_{y \in N(x)\cap V_j}C_{x,y}(\textcolor{black}{g}^{(n)}_{j}(x)-\textcolor{black}{g}^{(n)}_{j}(y))\vert}{m(x)}\leq 2$
for any $x \in V_{j}\setminus V_{n}$. 
\medskip

By applying {the} maximal principle to the function $h^{(n)}_{j}$, we observe $\max_{x \in V_{j}}\vert {h^{(n)}_{j}}(x)\vert \leq 2\sqrt{m(V)K}
({\rm deg^+}(D)+{\rm deg^-}(D)+1)/\min\{m(x)\mid x \in V_{n}\}$ which implies 
$\max_{x \in V_{j}}
\linebreak
\vert g^{(n)}_{j}(x)\vert \leq 2\sqrt{m(V)K}({\rm deg^+}(D)+{\rm deg^-}(D)+1)/\min\{m(x)\mid x \in V_{n}\}+1$.
\medskip

Accordingly, we see that the function $f^{(n)}_{j}$ taking $f_{n}$ on $V_{n}$ 
and $g^{(n)}_{j}$ on $V_{n}^c$ enjoys 
$\Vert \frac{1}{m}\Delta_{j} f^{(n)}_{j}1_{V_{n}^c}\Vert_{L^1(\nu)}
\leq 2m(V_{j}\setminus V_{n})$. In other words, $f^{(n)}_{j}$ meets the conditions (i)-(iii) in the assertion. \qed

\medskip
\medskip

{In what follows, }we discuss divisor $D=\sum_{x \in V_G} \ell(x)i(x)1_{\{x\}}$
satisfying $D=\sum_{x \in V_G} \vert \ell(x)\vert i(x)<\infty$.
In the next lemma, we start with a divisor $(D)_{n_0(\varepsilon)}$ satisfying ${\rm deg}^+(D-(D)_{n_0(\varepsilon)})+{\rm deg}^-(D-(D)_{n_0(\varepsilon)})< \varepsilon$ for a given $\varepsilon>0$
and we only focus on the subgraphs \textcolor{black}{$G_n=(V_n,E_n)$ with $n\geq n_0(\varepsilon)$}.
\textcolor{black}{Here and in the sequel, we assume that 
\Erase{$n_0(\varepsilon) \geq \hbox{\Erase{$\max\{n_0$,}} n_G
$ where $n_G$ appears at (\ref{n_G}) in the proof of Theorem \ref{poincare}.} 
\textcolor{black}{$n_0(\varepsilon)$ is sufficiently large as mentioned at the beginning of this section. }\Erase{such that $\rho_n<A$ holds for $n\ge n_0(\varepsilon)$ in Theorem \ref{poincare}.}} 
For any pair of positive integers $j,n$ with 
$j>n
\geq n_0(\varepsilon)$
and total order $O_{j} \in \mathcal O_{j}$, 
 the restriction of {the} total order of $O_{j}$ to $V_n$ is denoted by $O_{j}\vert_{V_n}$.
 \medskip

\textcolor{black}{\Erase{Now} Here, we make an attempt on choice of subsequence $\{n_k\}$
of positive integers so that the consistency in the sense $O_{n_j}\vert_{V_k}= O_{n_k}$ is satisfied for any pair $j,k$ of integers
with $j>k$.} 
%
\textcolor{black}{
Then, not only is the convergence of subsequence of $\{r_j((D)_n)\}$ for 
an arbitrarily fixed $n$ required,
but a procedure of taking limit as $j \to \infty$ and $n \to \infty$ in the doubly indexed sequence
$\{r_j((D)_n)\}$ should be reasonably organized to determine a characteristic value of such a divisor $D$ for our
Riemann-Roch theorem.  In the following lemma, we obtain such reasonable subsequence that we can eventually 
obtain the limit which makes sense for that purpose.
}
 
  \medskip
 \begin{Lem}\label{n_k}
 For any $\varepsilon >0$, there exists a sequence $\{O_{N(\varepsilon/2^{l})}\}$ of total orders 
 satisfying $O_{N(\varepsilon/2^{j})}\in \mathcal{O}_{N(\varepsilon/2^{j})}$ with $m(V_{N(\varepsilon/2^{j})}^c)<\varepsilon/2^{j}$
 for any non-negative integer $j$ and a sequence $\{n_{j}\}$ satisfying $n_1< n_2 < \dots$ and $n_{j+1} \geq  N(\varepsilon/2^{j})$ 
for any non-negative integer $j$
 such that  
\begin{equation}\label{ordeq}
r_{n_k}((D)_{n_{l}})
=\Big( \min_{D'\overset {n_k}{ \sim} (D)_{n_{l}}, O_{N(\varepsilon/2^{l})} =O_{n_k}\vert_{V_{N(\varepsilon/2^{l})}}, O_{n_k} \in \mathcal O_{n_k}}
{\rm deg}^{+}(D' - \nu_{O_{n_k}})\Big)-i_{(G_{n_k},C_{n_k})},
\end{equation}
whenever $k>l$. In particular,  $k>l$ implies
$O_{N(\varepsilon/2^{l})} = O_{n_k}\vert_{V_{N(\varepsilon/2^{l})} }$ and
\begin{equation}\label{deg2}
{\rm deg}^+(\nu_{O_{n_{l}}}-\nu_{O_{n_k}})+{\rm deg}^-(\nu_{O_{n_{l}}}-\nu_{O_{n_k}}) <m(V_{N(\varepsilon/2^{\min\{k,l\}})}^c))<\varepsilon/2^{\min\{k,l\}}
\end{equation}
\noindent
for any positive integers $k$ and $l$.
 \end{Lem}

 \medskip
{\it  Proof.} 
We take a divisor $(D)_{n_0(\varepsilon)}$ satisfying ${\rm deg}^+(D-(D)_{n_0(\varepsilon)})+{\rm deg}^-(D-(D)_{n_0(\varepsilon)})< \varepsilon$.
Since the finiteness $m(V)<\infty$ implies $\lim_{n\to \infty} m(V_n)=m(V)$,
there exists a positive integer $N(\varepsilon)$ such that $m(V_{N(\varepsilon)}^c)<\varepsilon$.
We may assume that $n_0(\varepsilon)\geq N(\varepsilon)$.
Since the cardinality of $\mathcal O_{N(\varepsilon)}$ is finite, the sequence $O_{N(\varepsilon)+1} \in \mathcal O_{N(\varepsilon)+1},O_{N(\varepsilon)+2} \in \mathcal O_{N(\varepsilon)+2},\dots$ admits
a subsequence $O_{n_1(\varepsilon)},O_{n_2(\varepsilon)},\dots$ with $N(\varepsilon)\leq n_0(\varepsilon)\leq n_1(\varepsilon)<n_2(\varepsilon)<\dots$ such that
$O_{N(\varepsilon)} =O_{n_k(\varepsilon)}\vert_{V_{N(\varepsilon)}}$ and 

\begin{align*}
r_{n_k(\varepsilon)}((D)_{n_0(\varepsilon)})
=&\Big( \min_{D'\overset {n_k(\varepsilon)}{ \sim} (D)_{n_0(\varepsilon)}, O_{N(\varepsilon)} =O_{n_k(\varepsilon)}\vert_{V_{N(\varepsilon)}}, O_{n_k(\varepsilon)} \in \mathcal O_{n_k(\varepsilon)}}
{\rm deg}^{+}(D' - \nu_{O_{n_k(\varepsilon)}})\Big)\\
&-i_{(G_{n_k(\varepsilon)},C_{n_k(\varepsilon)})}
\end{align*}

\noindent
for any $k \geq 1$.
We may assume that 
${\rm deg}^+(D-(D)_{n_1(\varepsilon)})+{\rm deg}^-(D-(D)_{n_1(\varepsilon)})< \varepsilon/2$

\medskip

By taking sufficiently large $n_2(\varepsilon)$, we may concentrate our attention to the case $n_2(\varepsilon)\geq N(\varepsilon/2)$.
Since the cardinality of $\mathcal O_{N(\varepsilon/2)}$ is finite, the sequence $O_{n_2(\varepsilon)},
O_{n_3(\varepsilon)}, 
\linebreak
O_{n_4(\varepsilon)},\dots$ admits
a subsequence $O_{n_2(\varepsilon/2)},O_{n_3(\varepsilon/2)},\dots$ with $N(\varepsilon/2)\leq n_2(\varepsilon/2)<n_3(\varepsilon/2)<\dots$ such that
$O_{N(\varepsilon/2)} =O_{n_k(\varepsilon/2)}\vert_{V_{N(\varepsilon/2)}}$ and

\begin{align*}
r_{n_k(\varepsilon/2)}((D)_{n_1(\varepsilon)})
=&\Big( \min_{D'\overset {n_k(\varepsilon/2)}{ \sim} (D)_{n_1(\varepsilon)}, O_{N(\varepsilon/2)} =O_{n_k(\varepsilon/2)}\vert_{V_{N(\varepsilon/2)}}, O_{n_k(\varepsilon/2)} \in \mathcal O_{n_k(\varepsilon/2)}}
{\rm deg}^{+}(D' - \nu_{O_{n_k(\varepsilon/2)}})\Big)\\
&-i_{(G_{n_k(\varepsilon/2)},C_{n_k(\varepsilon/2)})}
\end{align*}

\noindent
for any $k \geq 2$. 
We may assume that 
${\rm deg}^+(D-(D)_{n_2(\varepsilon/2)})+{\rm deg}^-(D-(D)_{n_2(\varepsilon/2)})< \varepsilon/4$

\medskip

By repeating this procedure, we obtain a subsequence $n_{j+1}(\varepsilon/2^{j}),n_{j+2}(\varepsilon/2^{j}),\cdots$ of 
$n_{j+1}(\varepsilon/2^{j-1}),n_{j+2}(\varepsilon/2^{j-1}),\cdots$
with $N(\varepsilon/2^{j})\leq n_{j+1}(\varepsilon/2^{j})<n_{j+2}(\varepsilon/2^{j})<\dots$ such that
$O_{N(\varepsilon/2^{j})} =O_{n_k(\varepsilon/2^{j})}\vert_{V_{N(\varepsilon/2^{j})}}$ and 

\begin{align*}
&r_{n_k(\varepsilon/2^{j})}((D)_{n_{j}(\varepsilon/2^{j-1})})\\
&=\Big( \min_{D'\overset {n_k(\varepsilon/2^{j})}{ \sim} (D)_{n_{j}(\varepsilon/2^{j-1})}, O_{N(\varepsilon/2^{j})} =O_{n_k(\varepsilon/2^{j})}\vert_{V_{N(\varepsilon/2^{j})}}, O_{n_k(\varepsilon/2^{j})} \in \mathcal O_{n_k(\varepsilon/2^{j})}}
{\rm deg}^{+}(D' - \nu_{O_{n_k(\varepsilon/2^{j})}})\Big)\\
&\quad -i_{(G_{n_k(\varepsilon/2^{j})},C_{n_k(\varepsilon/2^{j})})}
\end{align*}

\noindent
for any $k > j$. 
We may assume that 
${\rm deg}^+(D-(D)_{n_{j+1}(\varepsilon/2^{j})})+{\rm deg}^-(D-(D)_{n_{j+1}(\varepsilon/2^{j})})
< \varepsilon/2^{j}$.
 \medskip
 
As a result, by taking $n_1=n_1(\varepsilon),n_2=n_2(\varepsilon/2),\dots$, we obtain the sequence $n_1,n_2,\dots$ which meets all conditions 
in the assertion.

 \medskip
 
(\ref{deg2}) follows from the straightforward estimate $\vert \nu_{O_{n}}(x)-\nu_{O_{n'}}(x)\vert \leq m(x)$ for any $x \in V_n\cap V_{n'}$ with
$O_{n}\vert_{V_{N(\varepsilon)}}=O_{n'}\vert_{V_{N(\varepsilon)}}$. \qed

%
%
%
%
%
  \medskip
  
\textcolor{black}{We are now in position to assert the convergence of the sequence $\{r_{n_k}((D)_{n_l})\}_{k = l+1}^{\infty}$
as $k \to \infty$, for any fixed positive integer $l$.}

\medskip

\medskip
 \begin{Prop}\label{convergence} 
If 
\begin{equation}\label{ordeq2} 
\rho_{n}
m(S_{n})/\min_{x \in V_{n}}m(x)\to 0 \mbox{ as }n\to 
\infty,
\end{equation}
 then $\{r_{n_k}((D)_{n_l})\}_{k = l+1}^{\infty}$ converges
as $k \to \infty$, for any fixed non-negative integer $l$, where $n_1,n_2,\dots$ is the subsequence satisfying (\ref{ordeq})
 associated with a sequence $\{O_{N(\varepsilon/2^{j})}\}$ of total orders in Lemma \ref{n_k}. 
\end{Prop}

\medskip

{\it Proof.}   \textcolor{black}{First we note that, the assumptions in Theorem  \ref{poincare} and Lemma \ref{spectral gap} are satisfied since (\ref{ordeq2})} implies $\limsup_{n\to \infty}\rho_n< \min\{A,\frac{\lambda_G}{\lambda_G+1}\}$.
\Erase{Let  $\{f_{j}\}_{j=n}^{\infty}$ be the sequence consisting of $\mathbb Z$-valued functions 
each of which is  the minimizer for $r_j(D)$ taken in Lemma \ref{subsequence} with $j\geq n$.}
\Erase{We take the sequence $n_1,n_2,\dots$ satisfying (\ref{ordeq}) associated with the 
sequence $\{O_{N(\varepsilon/2^{l})}\}$ of total orders in Lemma \ref{n_k}} 
\medskip

\textcolor{black}{On the other hand, we also note that it suffices to show the assertion in the case $l=1$ because the following proof can be performed 
for any other positive integer $l$. In what follows, $(D)_{n_1}$ will be denoted simply by $D$.} \textcolor{black}{From now on, we use notation $\ell$ 
persistently for assigning an integer \textcolor{black}{$\ell$} satisfying  $ \ell > k$ so that confusion with the use of notation $l$ assigning integer less than $k$ \Erase{so far }is avoided\Erase{.
F} for comparing $r_{n_{\ell}}(D)$ with $r_{n_{k}}(D)$. Here, we take a minimizer $f_{n_{\ell}} $ for $r_{n_{\ell}}(D)$ and a minimizer $f_{n_{k}} $  for the other one.}
\medskip

By applying the identity just before
 (\ref{laplacian}) for the operators $\Delta_{n_{\ell}}$ and $\Delta_{n_{k}}$ instead of $L$ and $L_U$, we have
 $\Delta_{n_{\ell}}\phi(x)= \Delta_{n_{k}} \phi(x) + \big(\sum_{y\in N(x)\cap S_{n_k+1}}C_{x,y}(\phi(x)-\phi(y))\big)1_{S_{n_{k}}}
 +\Delta_{n_{\ell}}\phi(x)1_{V_{n_{\ell}}\setminus V_{n_{k}}}$ on $V_{n_k}$ for any function $\phi$ defined on $V_{n_{\ell}}$, 
where we recall
$S_{n_k}=\{x | d(v_0,x)=n_k\}$. 
 Accordingly, for the 
 function $f^{(n_k)}_{n_{\ell}}$ taken in Lemma \ref{(I)} by starting with the minimizer $f_{n_k}$, we have
\begin{align*}
r_{n_{\ell}}(D) &= {\rm deg}^+(D_{n_{\ell}}'-\nu_{O_{n_{\ell}}}) - i(G_{n_{\ell}},C_{n_{\ell}})\\
&={\rm deg}^+(D+\Delta_{n_{\ell}} f_{n_{\ell}}-\nu_{O_{n_{\ell}}})- i(G_{n_{\ell}},C_{n_{\ell}})\\
&\le {\rm deg}^+(D+\Delta_{n_{\ell}} f^{(n_k)}_{n_{\ell}}-\nu_{O_{n_{\ell}}}) - i(G_{n_k},C_{n_k})+\vert  i(G_{n_{\ell}},C_{n_{\ell}})-i(G_{n_k},C_{n_k})\vert\\
&= {\rm deg}^+(D+\Delta_{n_{k}} f_{n_{k}}+ \Big(\sum_{y\in N(x)\cap S_{n_k+1}}C_{x,y}( f^{(n_k)}_{n_{\ell}}(x)- f^{(n_k)}_{n_{\ell}}(y))\Big)1_{S_{n_{k}}}\\
&\quad +\Delta_{n_{\ell}}  f^{(n_k)}_{n_{\ell}}1_{V_{n_k}^c}-\nu_{O_{n_{\ell}}}) - i(G_{n_k},C_{n_k})+\vert  i(G_{n_{\ell}},C_{n_{\ell}})-i(G_{n_k},C_{n_k})\vert.
\end{align*}
Lemma \ref{(I)} shows that

\begin{align*}
&{\rm deg}^{+}(\sum_{x \in S_{n_k}}\sum_{y\in N(x)\cap S_{n_k+1}}C_{x,y}(f^{(n_k)}_{n_{\ell}}(x)-f^{(n_k)}_{n_{\ell}}(y)))\\ 
&+{\rm deg}^{-}(\sum_{x \in S_{n_k}}\sum_{y\in N(x)\cap S_{n_k+1}}C_{x,y}(f^{(n_k)}_{n_{\ell}}(x)-f^{(n_k)}_{n_{\ell}}(y))) \\ 
&\le \sum_{x\in S_{n_k}}\sum_{y\in N(x)\cap S_{n_k+1}} C_{x,y}(\vert f^{(n_k)}_{n_{\ell}}(x)\vert+\vert f^{(n_k)}_{n_{\ell}}(y)\vert)  \\ 
&\le 4\big(\sqrt{m(V)K}({\rm deg^+}(D)+{\rm deg^-}(D)+1)/\min\{m(x)\mid x \in V_{n_k}\}+1\big)\\
&\qquad \times\sum_{x\in S_{n_k}}\sum_{y\in N(x)\cap S_{n_k+1}}C_{x,y} \\
&\le  \rho_{n_k}m(S_{n_{k}})
\times 4\big((\sqrt{m(V)K}({\rm deg^+}(D)+{\rm deg^-}(D)+1)/\min\{m(x)\mid x \in V_{n_k}\}+1\big).
\end{align*}

Since our estimate on $\Delta_{j} f^{(n)}_{j}$ obtained in (iii) of Lemma \ref{(I)} is valid for $\Delta_{n_{\ell}} f^{(n_k)}_{n_{\ell}}$,
we see

$$\lim_{k\to \infty}{\sup_{\ell > k}\rm deg}^{\pm}\Delta_{n_{\ell}} f^{(n_k)}_{n_{\ell}}1_{V_{n_{\ell}}\setminus V_{n_k}}=0.$$

 \noindent
Hence, 
\begin{align*}
r_{n_{\ell}}(D)\leq & r_{n_k}(D)+{\rm deg}^+(\nu_{O_{n_{\ell}}}-\nu_{O_{n_k}})+{\rm deg}^-(\nu_{O_{n_{\ell}}}-\nu_{O_{n_k}})
+ \vert i_{(G_{n_{\ell}},C_{n_{\ell}})} - i_{(G_{n_k},C_{n_k})}\vert \\
&+\rho_{n_k}m(S_{n_{k}} )
4\big(\sqrt{m(V)K}({\rm deg^+}(D)+{\rm deg^-}(D)+1)/\min\{m(x)\mid x \in V_{n_k}\}+1\big)\\
&+o(1).
\end{align*}
Combining this and (\ref{ordeq2}) with ${\rm deg}^+(\nu_{O_{n_{\ell}}}-\nu_{O_{n_k}})+{\rm deg}^-(\nu_{O_{n_{\ell}}}-\nu_{O_{n_k}})<\varepsilon/2^{k}$ as obtained in (\ref{deg2}) and $\vert i_{(G_{n_{\ell}},C_{n_{\ell}})} - i_{(G_{n_k},C_{n_k})}\vert\to 0\ (\mbox{as }k,\ell\to\infty)$,
it turns out that $r_{n_{\ell}}(D)\le r_{n_k}(D)+\varepsilon/2^{k} + o(1)$ as $\ell\to\infty$ for any $k$,
which implies $\limsup_{\ell \to \infty}r_{n_{\ell}}(D)\leq  \liminf_{k\to \infty}r_{n_k}(D)$, in other words, $r_{n_k}(D)$ converges as $k\to\infty$.\qed

\medskip
The limit in this  proposition depends on the choice of the sequence $\{O_{n_k}\}$ of total orders. However,  
 as long as a divisor $D$ is supported by a finite graph, we can define
\[
r_{\{O_{n_k}\}}(D)=\lim_{k\to\infty} r_{n_k}(D)
\]
\noindent
by taking a subsequence $\{ O_{n_k}\}$ of total orders  as in Lemma \ref{n_k}.

\medskip

\begin{Remark}
If one takes another 
base \Erase{{\rm (}reference{\rm )}}
vertex $v'_0$ satisfying the condition (\ref{ordeq2}) in Proposition \ref{convergence}, instead of $v_0$, 
then $r'_{n_{k'}}(D)$ is defined \Erase{for}{as} the 
same
divisor $D$ 
\Erase{as}{in} 
the proposition. By taking a similar procedure in the proof 
of the proposition, one sees not only $\limsup_{\textcolor{black}{\ell'} \to \infty}r'_{n_{\ell'}}(D)\leq  \liminf_{k\to \infty}r_{n_k}(D)$
but {also} $\limsup_{\ell \to \infty}r_{n_{\ell}}(D)\leq  \liminf_{\textcolor{black}{k'}\to \infty}r'_{n_{k'}}(D)$. Accordingly, $\lim_{k\to\infty} r_{n_k}(D)$
does not depend on the choice of the 
base \Erase{{\rm (}reference{\rm )}}
vertex satisfying (\ref{ordeq2}).
\end{Remark}

\medskip
 For a divisor $D=\sum_{x \in V_G} \ell(x)i(x)1_{\{x\}}$ 
 on $G$,
we introduce effective divisors $D^+$ and $D^-$ given respectively by $D^+=\sum_{x \in V_G, \ell(x)>0} \ell(x)i(x)1_{\{x\}}$ and $D^-=-\sum_{x \in V_G, \ell(x)<0} \linebreak\ell(x)i(x)1_{\{x\}}$.  We
recall that its restriction $\sum_{x \in V_{n-1}} \ell(x)i(x)1_{\{x\}}$ to $V_n$ is denoted by $(D)_n$. 
\medskip

\begin{Cor} For any divisor $D=\sum_{x \in V_G} \ell(x)i(x)1_{\{x\}}$ on $G$ 
satisfying $\sum_{x \in 
V_G
}\vert \ell(x)\vert \linebreak i(x) < \infty$, $\{r_{\{O_{n_k}\}}((D)_{n_{l}})\}$ is a Cauchy sequence, where $\{(D)_{n_{l}}\}$ stands for a sequence of divisors in Lemma \ref{n_k}.

\end{Cor}

{\it Proof.}
By applying Lemma \ref{fundamental} to effective divisors $((D)_{n_{l'}}-(D)_{n_{l}})^{\pm}$ with $n_{l'}>n_{l}$, we have
$$r_{n_k}((D)_{n_{l}}) \leq r_{n_k}((D)_{n_{l}}+((D)_{n_{l'}}-(D)_{n_{l}})^{+})\leq r_{n_k}((D)_{n_{l}})+{\rm deg}^{+}((D)_{n_{l'}}-(D)_{n_{l}})$$

\noindent
and 
$$r_{n_k}((D)_{n_{l}})-{\rm deg}^{-}((D)_{n_{l'}}-(D)_{n_{l}}) \leq r_{n_k}((D)_{n_{l'}})\leq r_{n_k}((D)_{n_{l}})+{\rm deg}^{+}((D)_{n_{l'}}-(D)_{n_{l}}).$$

\medskip
\textcolor{black}{Thanks to Proposition \ref{convergence}, }
by passing the limit as $k\to \infty$, we have 
$$\vert r_{\{O_{n_k}\}}((D)_{n_{l}})- r_{\{O_{n_k}\}}((D)_{n_{l'}})\vert \leq {\rm deg}^{+}((D)_{n_{l'}}-(D)_{n_{l}})+{\rm deg}^{-}((D)_{n_{l'}}-(D)_{n_{l}}).$$

\medskip
Thanks to the finiteness $\sum_{x \in V_{G}}\vert \ell(x)\vert i(x) < \infty$,  for any $\varepsilon>0$,
there exists a positive integer $n_0(\varepsilon)$ such that 
$\sum_{x \notin V_{n_0(\varepsilon)}}\vert \ell(x)\vert i(x)={\rm deg}^+(D-(D)_{n_0(\varepsilon)})+{\rm deg}^-(D-(D)_{n_0(\varepsilon)})<\varepsilon$. Accordingly, it turns out that
$$l'>l\geq n_0(\varepsilon) \Rightarrow \vert r_{\{O_{n_k}\}}((D)_{n_{l}})-r_{\{O_{n_k}\}}((D)_{n_{l'}}))\vert <  \varepsilon$$

\noindent
independently of the choice of the sequence $\{O_{n_k}\}$. \qed

\medskip

This  implies the convergence of the sequence $\{r_{\{O_{n_k}\}}((D)_{n_{l}})\}_{l=1}^{\infty}$
and allows us to define $r(D)=\inf_{\{O_{n_k}\}}\lim_{l\to\infty}$ $r_{\{O_{n_k}\}}((D)_{n_{l}})$
for any divisor $D$ satisfying  $\sum_{x \in V_G}\vert \ell(x)\vert \linebreak i(x)<\infty$. 
Thanks to the finiteness of the measure $m$, by the definition of the characteristic  $\mathfrak e_{(G_n,C_n)}$ 
in the second section, it is easy to see that the sequence $\{\mathfrak e_{(G_n,C_n)}\}$ converges 
and the limit is independent of choice of the base \textcolor{black}{\Erase{(reference) }}vertex $v_0$.
\medskip

\begin{Thm} 
Let $G=(V_G,E_G)$ be a locally finite connected graph of finite volume satisfying (\ref{ordeq2}). 
For any divisor $D$ with 
${\rm deg}^{+}(D)+{\rm deg}^{-}(D)<\infty$, 
the Riemann Roch theorem holds on $G$: 
$$r(D)-r(K_G-D) = {\rm deg}(D) + \mathfrak e_{(
G
,C)},$$
where
$\mathfrak e_{(G,C)}=\lim_{n\to\infty}\mathfrak e_{(G_n,C_n)}$.
\end{Thm}

{\it Proof.} For a divisor $D$ with ${\rm deg}^{+}(D)+{\rm deg}^{-}(D)<\infty$, we take a sequence $\{(D)_{n_{l}}\}$
associated with $D$ in the sense of Lemma \ref{n_k}.
Riemann-Roch theorem
on finite weighted graph{s} shows that $k \geq l$ implies
$$
r_{n_k}((D)_{n_{l}})-r_{n_k}(K_{G_{n_k}}-(D)_{n_{l}}) = {\rm deg}((D)_{n_{l}}) + \mathfrak e_{(
G_{n_k}
,C_{n_k})}.
$$
From Corollary 2.4 and Proposition 4.4 starting with $(D)_{n_{l}}$ \Erase{instead of $D$}, by
letting $k\to\infty$, we can derive
$$
r_{\{O_{n_k}\}}((D)_{n_{l}})-r_{\{O_{n_k}\}}(K_{G}-(D)_{n_{l}}) = {\rm deg}((D)_{n_{l}}) + \mathfrak e_{(G,C)}.
$$

\noindent
By passing the limit as $l\to \infty$ and taking the infimum over all \textcolor{black}{\Erase{choice}} of sequence{s} $\{O_{n_k}\}$ of total orders
in Lemma \ref{n_k}, it turns out that

$$r(D)-r(K_G-D) = {\rm deg}(D) + \mathfrak e_{(G,C)}.$$

\qed

\label{}





\end{document}